\documentclass[12pt]{article}
\usepackage{a4,amssymb}
\usepackage{graphicx}
\def\Bbb{\mathbb}

\title{\bf Semistable K3-surfaces \\ with icosahedral symmetry}
\author{Jan Stevens}
\date{}

\def\wtx{\widetilde{X}}
\def\wl#1{\overline{#1}}
\def\kop#1\par{\par\bigskip{\bf #1\unskip.}\par\medskip}
\def\pab#1{{\partial\over \partial #1}}

\def\ct{{\cal T}}
\def\cg{{\cal G}}
\def\cf{{\cal F}}
\def\al{\alpha}
\def\be{\beta}
\def\la{\lambda}
\def\ga{\gamma}
\def\de{\delta}
\def\De{\Delta}
\def\Dee{\Delta}
\def\ep{\varepsilon}
\def\th{\vartheta}

\def\sier#1{{\cal O}_{#1}}
\def\breuk#1/#2{{\textstyle {#1\over #2}}}
\def\supp{\mathop{\rm supp}}
\def\cn{\mathop{\rm :}}
\def\half{\breuk1/2}
\def\C{{\Bbb C}}  \def\P{{\Bbb P}}  \def\Z{{\Bbb Z}} \def\N{{\Bbb N}}
\def\D{{\cal D}{\it iff}}
\def\coker{\mathop{\rm coker}}
\def\lra{\longrightarrow}
\def\nroep #1.{\protect\refstepcounter{subsection}\protect%
\medbreak\noindent{\bf (\thesubsection)} {\sl #1\/}.\enspace}
\def\roep #1.{\medbreak\noindent{\sl #1\/}.\enspace}
\def\endroep{\par\medbreak}
\def\thesection{\arabic{section}.}
\def\thesubsection{\thesection\arabic{subsection}}

\long\def\comment#1\endcomment{}
\newcommand{\neu}{\protect\refstepcounter{subsection}\protect\medbreak
\noindent{\bf (\thesubsection)} \enspace\ignorespaces}
\newcommand{\bfneu}[1]{\protect\refstepcounter{subsection}\protect\medbreak
\noindent{\bf (\thesubsection) #1} \enspace\ignorespaces}

\long\def\proclaim#1. #2\par{\protect\refstepcounter{subsection}\protect
\medbreak\noindent{\bf(\thesubsection) #1.}\enspace{\sl#2}\par
\ifdim\lastskip<\medskipamount \removelastskip\penalty55\medskip\fi}
\def\Box{\square}
\newcommand{\qed}{{
\unskip\nobreak\hfil\penalty50\hskip2em\hbox{}\nobreak\hfil$\Box$
\parfillskip=0pt \finalhyphendemerits=0 \par\bigbreak}}

\catcode`\@=11
\def\eqalign#1{\null\,\vcenter{\openup1\jot \m@th
  \ialign{\strut\hfil$\displaystyle{##}$&$\displaystyle{{}##}$\hfil
      \crcr#1\crcr}}\,}
\catcode`\@=12

\begin{document}
\maketitle
\begin{abstract}
\noindent In a Type III degeneration of $K3$-surfaces the dual graph of 
the central fibre is a triangulation of $S^2$. 
We realise the  tetrahedral, octahedral and especially the icosahedral
triangulation in  families of $K3$-surfaces,  preferably with the
associated symmetry groups acting. 
\end{abstract}
\comment
In a Type III degeneration of K3-surfaces the dual graph of 
the central fibre is a triangulation of the 2-sphere.
We realise the tetrahedral, octahedral and especially the icosahedral
triangulation in  families of K3-surfaces,  preferably with the
associated symmetry groups acting.
\endcomment

\section*{Introduction}
A degeneration of surfaces is a 1-parameter family with general fibre
a smooth complex surface.
The case of $K3$-surfaces has attracted a lot of attention. 
A nice discussion  is contained in the introductory
first paper \cite{FM} of the bundle \cite{SAGS}.
One usually allows base change and modifications to obtain good models.
After a ramified cover of the base and resolution of singularities
we can assume that the degeneration $f\colon {\cal X}\to S \ni 0$
is {\sl semistable\/}: the zero fibre  $X=f^{-1}(0)$ is a  reduced
divisor with (simple) normal crossings in the smooth manifold $\cal X$. 
Further modifications of a $K3$-degeneration lead to a minimal model,
which falls into one of three types.

In a Type III degeneration of $K3$-surfaces  the dual graph of 
the central fibre is a triangulation of $S^2$. In this paper I
construct an example with my favourite triangulation, the icosahedral one.
A substantial part is taken up by the
tetrahedral case, which is easier to handle and allows more explicit
results. A second purpose of this paper is to link general theory 
with concrete computations.

There are two obvious ways to realise a semistable degeneration with
prescribed combinatorial type.
The first is to start with a singular total space, but correct 
central fibre. An example is the coordinate tetrahedron 
$T=x_0x_1x_2x_3$ in $\P^3$, 
where we get a degeneration $f\colon {\cal X}\to \P^1 \ni 0$
by blowing up the base locus of the pencil spanned by
$T$ and a generic quartic. The total space has then 
24 $A_1$-singularities. This is a minimal model in the Mori category,
but one has to take a small resolution to get a smooth total space.
We can arrange that each plane of the tetrahedron is blown up 
in $6$ points. The dual graph of the central fibre remains the same.
 
The second method is to try to smooth the putative central fibre.
For the tetrahedron this can be done directly.
We glue together
four cubic surfaces  along triangles. This 
normal crossings variety satisfies the topological conditions to be
a central fibre (the triple point formula), but one needs also 
a more subtle analytic condition ($d$-semistability), which translates
into equations on the coefficients in the equations.
The necessary deformation theory in general has been developed by
Friedman \cite{Fr}. The central result is that smoothing is
always possible  in the $K3$-case.
This holds both for abstract and embedded deformations.
One typically obtains different degenerations from the two
constructions, which fill up 19 dimensional families in a 20 dimensional 
deformation space of the normal crossings $K3$. 

For the octahedron both methods can again be applied. We find
the correct space as anticanonical divisor in the toric threefold
given by a cube.
In his monograph Ulf Persson invites the reader `to find a degeneration
into a dodekahedron  of rational surfaces' \cite[p.~126]{Pe}.
For a construction according to the first method the double curve
must be an anti-canonical divisor on each component. The most natural choice
for a rational surface is then a Del Pezzo of degree five. This makes
that the dodecahedron itself has degree 60, and it is exactly 
such a dodecahedron, obtained by glueing twelve Del Pezzo surfaces, which
the second method smoothes.
Unfortunately the computations are too difficult to give explicit
formulas. The same holds for a related problem in less variables, 
smoothing the Stanley-Reisner ring of the icosahedron. A semistable model
of such a degeneration has as central fibre
a complexified football. The extra components come from singularities
of the total space.

The last example suggests that one can get a dodecahedron out of
a central fibre with less than 12 components. This requires a breaking
of the symmetry. By combining the first and the second method
I obtain in
(\ref{icodual}) an  explicit degeneration, 
whose general fibre is a smooth $K3$-surface of degree 12
in $\P^7$, with special fibre
consisting of 6 planes with triangles as double curves
and 3 quadric surfaces with rectangular double curve. Its total space
has three singularities, which are isomorphic to cones over Del Pezzo
surfaces of degree 5, and 18 $A_1$ points. 
The dual graph of the central fibre on a suitable smooth model
is the icosahedron.

This paper is organised as follows.
In the first Section I recall the results on degenerations
of $K$3-surfaces, in particular that one can always realise
a particularly nice model, the $(-1)$-form. Section 2
brings as illustration  detailed computations for tetrahedra.
The results fit in with the general deformation theory, which is
treated in the third Section, with special emphasis on degenerations
in $(-1)$-form. A short fourth Section introduces the combinatorial
tools to handle large systems of equations: the definitions
of Stanley-Reisner rings and Hodge algebras are reviewed.
The final section contains the dodecahedral degenerations.

\section{Semistable degenerations of $K$3-surfaces}
\neu
The name $K3$ has been explained by Andr\'e Weil: `en
l'honneur de Kummer, K\"ahler, Kodaira et de la belle montagne $K2$
au Cachemire' \cite[p.~546]{W}. He calls
any surface a $K3$, if it has the differentiable structure of a smooth 
quartic surface in $\P^3(\C)$. 
A Kummer surface is a quartic with 16 $A_1$-singularities. As these
singularities admit simultaneous resolution, the minimal resolution
of a Kummer surface deforms into a smooth quartic and is therefore
a $K3$ surface.
A quartic surface $X$ is simply connected,
so in particular $b_1(X)=0$ and has trivial canonical sheaf
by the adjunction formula: $X$ is an anti-canonical divisor in $\P^3$.
The modern definition of a $K3$-surface: $b_1(X)=0$ and $ K_X =0$, is 
equivalent with Weil's definition  because all
$K3$-surfaces form one connected family.
\neu
Let $f\colon {\cal X}\to S \ni 0$ be a proper surjective holomorphic
map of a 3-dimensional complex manifold $\cal X$ to a (germ of a) curve
$S$ such that the zero fibre  $X=f^{-1}(0)$ is a  reduced
divisor with (simple) normal crossings; then the degeneration $f$ is
called  {\sl semistable\/}. 

In the  $K3$ case the following holds (see \cite{FM}) for exact references):

\proclaim Theorem {\rm (Kulikov)}. Let $f\colon {\cal X}\to S$ be a
semistable degeneration of $K3$-surfaces. 
If all components of $X=f^{-1}(0)$ are K\"ahler, then there 
exists a modification ${\cal X}'$ of ${\cal X}$ such that
$K_{{\cal X}'}\equiv 0$.

A degeneration as in the conclusion of the theorem $(K_{{\cal X}}\equiv 0)$
is called a {\sl Kulikov model\/}.

\proclaim Theorem {\rm (Persson, Kulikov)}. Let $f\colon {\cal X}\to S$ 
be a Kulikov model of a degeneration of $K3$-surfaces with all
components of $X=f^{-1}(0)$  K\"ahler. 
Then either
\begin{itemize}
\item[\rm I] $X$ is smooth, or
\item[\rm II] $X$ is a chain of elliptic ruled components 
with rational surfaces at the ends and all double curves are smooth
elliptic curves, or
\item[\rm III] $X$ consists of rational surfaces meeting along rational curves
which form cycles on each component. The dual graph is a triangulation of
$S^2$.
\end{itemize}

According to the case division in the theorem one speaks of degenerations
of type  I, II, or III.
Without the K\"ahler assumption it is not always possible to arrange that
$K_{{\cal X}}\equiv 0$ \cite{Kull,N}. In particular it is possible that
the central fibre contains surfaces of type ${\rm VII}_0$.  Even
under the assumption $K_{{\cal X}}\equiv 0$ the list becomes longer
(see \cite[Thm.~2.1]{N}). The case that the central fibre contains
an Inoue-Hirzebruch surface is 
relevant for the deformation of cusp singularities \cite{Lo}.

A Kulikov model is not unique. The central fibre
can be modified with flops.
If $C^{-}$ is a smooth rational curve in $X$ with self intersection $-1$,
lying in a component $X_i$ and
intersecting the double curve transversally in one point lying in $X_j$,
then after the flop the curve $C^{+}$ lies in $X_j$. This operation is
also called elementary modification of type I along $C^{-}$. 
An elementary modification of type II is a flop in a curve $C^{-}$,
which is a component of the double curve and has self intersection $-1$
on both components $X_i$, $X_j$ on which it lies. There are two
triple points on $C^{-}$ involving the components $X_k$ and $X_l$.
After the flop $C^+$ is a double curve lying in the components 
$X_k$ and $X_l$.
Note that we we might lose projectivity by using elementary transformations.
\begin{figure}[h]
\centering
\includegraphics[width=10cm]{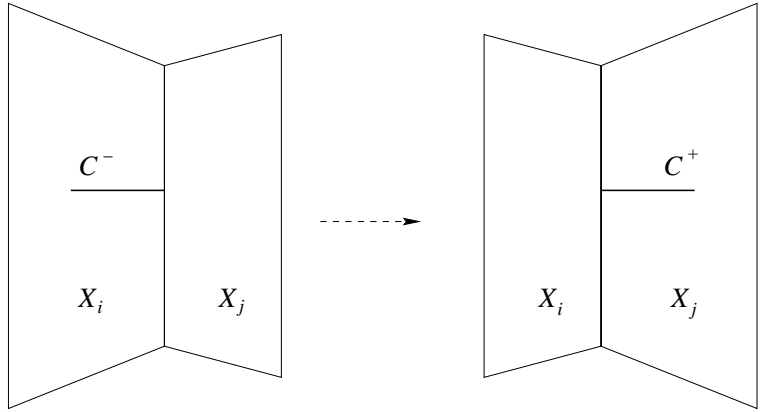}
\caption{Elementary modification of type I.}
\end{figure}
\begin{figure}[h]
\centering
\includegraphics[width=10cm]{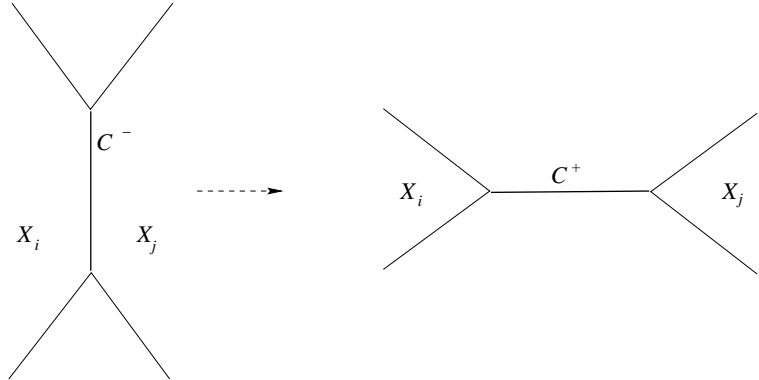}
\caption{Elementary modification of type II.}
\end{figure}

\proclaim The Minus One Theorem {\rm \cite{MM}}.
By modifications of type I and II one can achieve that 
every component of the double curve of the special fibre  
has self intersection $-1$ on both components on which it lies.

\neu
Let $X=\cup X_i$  be a normal crossings surface  with double locus
$D$. If $X$ is a divisor in a smooth 3-fold $M$ then one can define
the infinitesimal normal bundle $\sier D(X)$
as $\sier D(X)=\sier M(X)|_D$. It can be defined independent of $M$.
To this end,
let $I_{X_i}$ be the ideal sheaf of $X_i$ in $X$. It is locally
generated by one generator $z_i$, but not invertible as $z_i$ is
a zero divisor in $\sier X$. However $I_{X_i}|_D$ is locally free
\cite[(1.8)]{Fr}. Following Friedman one makes the following definitions.

\nroep Definition. The {\sl infinitesimal normal bundle\/} $\sier D(X)$ is the
line bundle dual to $\sier D(-X)$, where
$$
\sier D(-X)
= (I_{X_1}|_D) \otimes_{\sier D} \cdots \otimes_{\sier D} (I_{X_k}|_D) \;.
$$
\endroep

If $X$ is a divisor in $M$ the so defined bundle equals $\sier M(X)|_D$.
In particular, if $X$ is a central fibre
in a semistable degeneration ${\cal X}\to S$, then 
$\sier {\cal X}(X)\equiv \sier{\cal X}$ so $\sier D(X)=\sier D$. 
This gives a necessary condition for being a central fibre.

\nroep Definition. The normal crossings surface $X$ is {\sl $d$-semistable\/}
if $\sier D(X)=\sier D$.
\endroep

A consequence is the triple point formula:
let $D_{ij}=X_i\cap X_j$ and denote by $(D_{ij})^2_{X_i}$ the self 
intersection of $D_{ij}$ on $X_i$ and by $T_{ij}$ the number of
triple points on $D_{ij}$.
Then (cf.~\cite[Cor.~2.4.2]{Pe})
$$
(D_{ij})^2_{X_i}+(D_{ij})^2_{X_j}+T_{ij}=0\;.
$$

\nroep Definition.
A compact normal crossings surface is a {\sl $d$-semistable
$K3$-surface of type\/} III if $X$ is $d$-semistable, $\omega_X=\sier X$
and each $X_i$ is rational, the double curves $D_i\subset X_i$ 
are cycles of rational curves and the dual graph triangulates $S^2$.
If the conclusions of the Minus One Theorem hold, 
that  every component of the double curve has self intersection
$-1$ on either component of $X$ on which it lies, the surface $X$
is said to be in {\sl $-1$ form\/}.
\endroep

\section{Tetrahedra}

\neu
To realise a tetrahedron we start out with four general planes in 3-space.
They do not form a $d$-semistable
$K3$-surface, but the dual graph is a tetrahedron.
To write down a degeneration with this
special fibre we just take the pencil spanned by $T=x_0x_1x_2x_3$
and a smooth quartic. The symmetry group of
the tetrahedron (including reflections) acts if
we only take $S_4$-invariant quartics:
$$
Q=a \sigma_1^4 + b \sigma_1^2\sigma_2 + c \sigma_2^2 + d\sigma_1\sigma_3\;,
$$
where the $\sigma_i$ are the elementary symmetric functions in the
four variables $x_i$ and $a$, $b$, $c$ and $d$ are constants.

To obtain a family $f\colon {\cal X}\to S$ one has to blow up 
the base locus of the pencil. This can be done in several ways.
Blowing up $T=Q=0$ gives a total space which is singular, with 
in general 24 ordinary double points coming from the 24 intersection points
of $Q$ with the double curve of the tetrahedron $T$.
Arguably, this is the nicest model, and the best one can hope for
in view of the theory of minimal models of 3-folds.
A smooth model 
is obtained by a suitable 
small resolution of the 24 singularities.
The quartic intersects each edge of the
tetrahedron in four points. To get the $(-1)$-form two of them
have to be blown up in one face and the other two in the other face.
The central fibre then consists of four Del Pezzo surfaces of degree $3$.

Alternatively one can blow up the irreducible components of $T=Q=0$
one at a time. The advantage is that one has a projective model.
However it is not in $(-1)$-form and furthermore the symmetry is
not preserved. To achieve $(-1)$-form we have to apply 
modifications of type I; here we may loose projectivity.

\bfneu{The tetrahedron of degree 12.}\label{tetra}
We glue together four Del Pezzo surfaces of degree $3$. Take coordinates
$x_1$, \dots, $x_4$, $y_1$, \dots, $y_4$ on $\P^7$. Let $\{i,j,k,l\}
=\{1,2,3,4\}$. The Del Pezzo surface $X_i$ lies in $y_i=x_j=x_k=x_l=0$
and has an equation of the form 
$$
y_jy_ky_l-x_iF_i(x_i,y_j,y_k,y_l)=0\;,
$$
where $F_i$ is a quadratic form; more specifically, 
$$
F_i= \sum_{\al\neq i} f_i^{\al\al}y_\al^2 
+\sum_{\al,\be\neq i} f_i^{\al\be}y_\al y_\be +
\sum_{\al\neq i} g_i^{\al}x_iy_\al+h_ix_i^2\;.
$$
The condition that $X_i$ is nonsingular
in the vertices of the triangle $x_i= y_jy_ky_l=0$ is that the
coefficients $f_i^{\al\al}$ in $F_i$ do not
vanish.  

The ideal of the tetrahedron $X=\bigcup_i X_i$ has 14 generators,
the 4 cubic Del Pezzo equations  $y_jy_ky_l-x_iF_i$ and
10 quadratic monomials: the six products $x_ix_j$ and 
the four products $x_iy_i$.
The relations among them are:
\begin{equation}
\label{tetrarel}
\eqalign{
(x_ix_j)x_k&-(x_ix_k)x_j    \cr
(x_iy_i)x_j&-(x_ix_j)y_i  \cr
(y_jy_ky_l-x_iF_i)x_l &- (x_ly_l)y_jy_k+(x_ix_l)F_i  \cr
(y_iy_jy_k-x_lF_l)y_l&-(y_jy_ky_l-x_iF_i)y_i -(x_ly_l)F_l+(x_iy_i)F_i  \;.
}
\end{equation}

\proclaim Proposition.
The tetrahedron $X$ is $d$-semistable if and only if the four equations
$$
f_k^{jj}f_l^{kk}f_j^{ll}-f_j^{kk}f_k^{ll}f_l^{jj}=0
$$
are satisfied.

\roep Proof.
We look at the chart $y_4=1$. Then $x_4=0$ and we have the equations
$x_ix_j$,  $x_iy_i$, $y_1y_2y_3$ and $y_iy_j-x_kF_k$. In all points 
near the origin
$y_1y_2y_3\mapsto 1$ is generator  of the infinitesimal normal bundle
$\sier D(X)$. We now look on the $y_2$-axis.
The equation $y_2y_1-x_3F_3$ shows that the section  $y_1y_2y_3\mapsto 1$
has a pole in the zeroes of $F_3$ restricted to the $y_2$-axis, and likewise
in the zeroes of $F_1$ (using the equation $y_2y_3-x_1F_1$).

This shows that the expression, given in homogeneous coordinates
by 
$$
{(f_1^{22}y_2^2+f_1^{24}y_2y_4+f_1^{44}y_4^2)
(f_3^{22}y_2^2+f_3^{24}y_2y_4+f_3^{44}y_4^2) \over
y_1y_2y_3y_4}\;,
$$
represents a nonvanishing holomorphic section of $\sier D(X)$
on the whole line $y_1=y_3$. Similar expression can be found for the
other edges of the tetrahedron.
We get a global section if and only if we can find a quaternary form $f$
of degree $4$ which restricts to a multiple of the above denominator
for each line. 

For each face we find from the following Lemma 
the condition that the 12 points in the corresponding hyperplane
are cut out by a quartic. We obtain four equations, which in fact
are not independent: under our assumption that all $f^i_{jj}$ are
different from 0 one can derive one equation from the remaining three.
They give necessary and sufficient conditions for the existence
of the quaternary quartic.
\qed

\proclaim Lemma. \label{conditie}%
Consider $3n$ points $P_{i,\al}$ with $P_{i,1}$, \dots, $P_{i,n}$
smooth points of the triangle $x_1x_2x_3=0$, lying on the side $x_i=0$
and given by the binary form 
$B_i(x_j,x_k)=\sum_{m=0}^n b_{im}x_j^mx_k^{n-m}$, where $(i,j,k)$ is a
cyclic permutation of $(1,2,3)$. These points are cut out by a ternary
form of degree $n$ if and only if
$$
b_{10}b_{20}b_{30}=b_{1n}b_{2n}b_{3n}\;.
$$

\roep Proof.
Suppose $A(x_1,x_2,x_3)=\sum_{l+m+p=n}a_{lmp}x_1^lx_2^mx_3^p$ cuts out
the points. Then $A(0,x_2,x_3)$ is proportional to $B_1(x_2,x_3)$, 
so $(a_{0n0}\cn a_{00n})=(b_{10} \cn b_{1n})$. Likewise we have that 
$(a_{00n}\cn a_{n00})=(b_{20} \cn b_{2n})$ and
$(a_{n00}\cn a_{0n0})=(b_{30} \cn b_{3n})$.
Multiplying these ratios gives the condition.

Conversely, to find $A$ we may suppose that $b_{10}=b_{3n}=1=a_{0n0}$
(as no point lies at one of the vertices).
We put $A(0,x_2,x_3)=B_1(x_2,x_3)$,  $A(x_1,x_2,0)=B_3(x_1,x_2)$.
We also can take $b_{20}=a_{00n}$. As $b_{1n}=a_{00n}$
the condition gives now $b_{2n}=b_{30}=a_{n00}$ and we can set
$A(x_1,0,x_3)=B_2(x_1,x_3)$. The remaining monomials in $A$ are
divisible by $x_1x_2x_3$ and do not matter.
\qed

\nroep Remark. It is not surprising that only the extremal coefficients
$b_{i0}$, $b_{in}$ are involved, as they depend only on the product
of the coordinates of the points. Ignoring the other coefficients
we rename: $b_{i0}=:b_{jk}$, $b_{in}=:b_{kj}$. The condition
becomes $b_{jk}b_{ki}b_{ij}=b_{ji}b_{ik}b_{kj}$,
which is the form used in the proposition above.
\endroep

\bfneu{Infinitesimal deformations.}
We compute embedded deformations mo\-dulo coordinate transformations.
To this end we look at the equations as defining the affine cone
$C(X)$ over $X$. We follow the standard procedure
(see e.g.~\cite{St2}):
given equations $f_i$, satisfying relations $\sum f_ir_{ij}=0$, 
we have to lift the equations to $F_i=f_i+\ep f_i'$
and the relations to $R_{ij}=r_{ij}+r_{ij}'$, satisfying 
$\sum  F_iR_{ij}\equiv0 \pmod {\ep^2}$. This means that we have to find
$f_i'$ such that $\sum f_i'r_{ij}$ lies in the ideal generated
by the $f_i$. 
Using undetermined coefficients this is a finite dimensional problem
for each degree.
The deformations of $C(X)$ in degree $0$ give 
embedded deformations of $X$ in $\P^7$, while those
in degree $<0$ have an interpretation in terms of extensions of $X$:
they tells us of which varieties $X$ is a hyperplane section.
Our main interest lies in the degree 0 deformations, but as preparation 
we first compute those  of negative degrees.

\proclaim Proposition.
The dimension of $T^1_{C(X)}(-2)$ equals 4 and $\dim T^1_{C(X)}(-1)=16$.
In case $X$ is $d$-semistable $\dim T^1_{C(X)}(0)=22$, otherwise it is $21$.

\roep Proof.
Degree $-2$:
we perturb the quadratic equations with constants and the cubic equations
with linear terms. Write $x_ix_j+a_{ij}$. The first type of the relations
(\ref{tetrarel}) 
then gives $a_{ij}x_k-a_{ik}x_j=0\in \sier{C(X)}$, so $a_{ik}=0$.
Also the equation $x_iy_i$ are not perturbed.
Consider $ y_jy_ky_l-x_iF_i+\sum a_i^\al x_\al + \sum b_i^\al y_\al$.
The third relation gives $a_i^jx_j^2+\sum_{\al\neq j} b_i^\al x_jy_\al=0$
so we conclude that all coefficients  vanish, except the $a_i^i$,
which we may choose arbitrary. The last type
of relation is then also satisfied.

Degree $-1$:
consider the perturbations 
$$
x_ix_j+\sum a_{ij}^\al x_\al + \sum b_{ij}^\al y_\al\;.
$$
In the local ring we obtain the equation
$$
a_{ij}^k x_k^2  + \sum_{\al\neq k} b_{ij}^\al x_k y_\al -
a_{ik}^k x_j^2  + \sum_{\al\neq j} b_{ik}^\al x_j y_\al =0\;,
$$
from which we get $b_{ij}^\al=0$, $a_{ij}^k=a_{ij}^l=0$.
We now put 
$$
x_iy_i+\sum a_{ii}^\al x_\al + \sum b_{ii}^\al y_\al\;.
$$
We find
$$
a_{ii}^j x_j^2 + \sum_{\al\neq j} b_{ii}^\al x_jy_\al
-a_{ij}^j x_jy_i=0\;.
$$
We conclude $a_{ii}^j=0$, $b_{ii}^j=0$ for all $j\neq i$ and
finally $a_{ij}^j=b_{ii}^i$. In particular  $a_{ij}^j$ is independent
of $j$. We can use the coordinate transformation $x_i \mapsto
x_i-b_{ii}^i$ to get rid of the $a_{ij}^j$-term. So the equations
$x_ix_j$ are not perturbed at all. This means that $x_iy_j$ is only
perturbed with the term $a_{ii}^ix_i$, which can be made to vanish
by coordinate transformations in the $y_i$-variables.
As above we find that the only allowable pertubations of the cubic
equation $ y_jy_ky_l-x_iF_i$ are those divisible by $x_i$.
As we have used all coordinate transformations, all
monomials $x_i^2$, $x_iy_j$ can occur. This makes the dimension of
$T^1(-1)$ into $4\times 4$.

Degree $0$:
we proceed in the same way by first considering
the perturbations 
$$
x_ix_j+\sum a_{ij}^\al x_\al^2 + \sum b_{ij}^{\al\be} x_\al y_\be
+\sum c_{ij}^\al y_\al^2 + \sum d_{ij}^{\al\be} y_\al y_\be
\;.
$$
Multiplied with $x_k$ this gives the following terms in the local ring
$$
a_{ij}^k x_k^3 + \sum_{\be\neq k} b_{ij}^{k\be} x_k^2 y_\be
+\sum_{\al\neq k} c_{ij}^\al x_k y_\al^2 
+ \sum_{\al,\be\neq k} d_{ij}^{\al\be} x_k y_\al y_\be
\;.
$$
We conclude that all coefficients occurring here vanish.
In particular 
$a_{ij}^k=0$. Using the coordinate transformations
$x_j \mapsto x_j-a_{ij}^ix_i$ we may suppose that all $a_{ij}^\al$
vanish. 
We are left with 
$$
x_ix_j+ \sum b_{ij}^{i\be} x_i y_\be + \sum b_{ij}^{j\be} x_j y_\be
+ d_{ij}^{kl} y_k y_l
\;.
$$
With the perturbations 
$$
x_iy_i+\sum_{\al\neq i} a_{ii}^\al x_\al^2 + 
\sum_{\al\neq i} b_{ii}^{\al\be} x_\al y_\be
+\sum_{\al\neq i} c_{ii}^\al y_\al^2 + 
\sum_{\al,\be\neq i} d_{ii}^{\al\be} y_\al y_\be\;,
$$
where we used coordinate transformations $x_i\mapsto x_i-c_{ii}^iy_i$,
$x_i\mapsto x_i-d_{ii}^{ij}y_j$, $y_i\mapsto y_i-a_{ii}^ix_i$ 
and $y_i\mapsto y_i-b_{ii}^{ij}y_j$ to
remove some coefficients,
we now get (using the $j$th Del Pezzo equation)
$$
a_{ii}^j x_j^3 + \sum b_{ii}^{j\be} x_j^2 y_\be
+\sum_{\al\neq i} c_{ii}^\al x_j y_\al^2 
+ \sum_{\al,\be\neq i} d_{ii}^{\al\be} x_j  y_\al y_\be
-\sum b_{ij}^{j\be} x_j y_i y_\be
- d_{ij}^{kl} x_j F_j=0\;.
$$
Using the explicit expression for $F_j$ we obtain the equations
$$
\begin{array}{ccccccccc}
a_{ii}^j &=&d_{ij}^{kl}h_j,\quad &
b_{ii}^{j\be} &=&d_{ij}^{kl} g_j^{\be},\quad &
c_{ii}^\al   &=& d_{ij}^{kl}f_j^{\al\al}, \\
-b_{ij}^{ji}   &=& d_{ij}^{kl}f_j^{ii},\quad &
d_{ii}^{kl}&=& d_{ij}^{kl}f_j^{kl},\quad &
-b_{ij}^{j\be}&=& d_{ij}^{kl}f_j^{i\be}.
\end{array}
$$
We can determine all coefficients, but because 
$c_{ii}^\al$ does not depend on $j$, we get two equations for it
$$
d_{ij}^{kl}f_j^{ll}=c_{ii}^l=d_{ik}^{jl}f_k^{ll}\;.
$$
We view these as one linear equation for the unknowns $d_{ij}^{kl}$.
The coefficient matrix of the resulting  linear system is 
$$
\pmatrix
{
0 & f_k^{jj} & -f_l^{jj} & 0 & 0 & 0 \cr
-f_j^{kk} & 0 &f_l^{kk}  & 0 & 0 & 0 \cr
f_j^{ll} &-f_k^{ll}  & 0 & 0 & 0 & 0 \cr
f_i^{kk} & 0 & 0 & 0 & -f_l^{kk} & 0 \cr
-f_i^{ll} & 0 & 0 & f_k^{ll} & 0 & 0 \cr
0 & 0 & 0 & -f_k^{ii} & f_l^{ii} & 0 \cr
0 & f_i^{ll} & 0 & -f_j^{ll} & 0 & 0 \cr
0 & -f_i^{jj} & 0 & 0 & 0 & f_l^{jj} \cr
0 & 0 & 0 &f_j^{ii}  & 0 & -f_l^{ii} \cr
0 & 0 & f_i^{jj} & 0 & 0 & -f_k^{jj} \cr
0 & 0 & -f_i^{kk} & 0 & f_j^{kk} & 0 \cr
0 & 0 & 0 & 0 & -f_j^{ii} & f_k^{ii} \cr
}
$$
It has a nontrivial solution if all $6\times6$ minors vanish.
Among those are
$$
f_j^{kk}f_k^{ll}f_l^{jj}(f_k^{jj}f_l^{kk}f_j^{ll}-f_j^{kk}f_k^{ll}f_l^{jj})
$$
from which we obtain that the square of 
$$
f_k^{jj}f_l^{kk}f_j^{ll}-f_j^{kk}f_k^{ll}f_l^{jj}
$$
lies in the ideal of the minors.
This is one of the four conditions for $d$-semistability.
There are three more  equations 
$$
f_k^{jj}f_j^{ll}f_l^{ii}f_i^{kk}-f_j^{kk}f_l^{jj}f_i^{ll}f_k^{ii}
$$
in the reduction of the ideal of minors,
which do not give new conditions if the $f_i^{jj}\neq0$, as
$$
\displaylines{
\qquad
f_l^{kk}(f_k^{jj}f_j^{ll}f_l^{ii}f_i^{kk}-f_j^{kk}f_l^{jj}f_i^{ll}f_k^{ii})
\hfill\cr\hfill{}=
f_l^{ii}f_i^{kk}(f_k^{jj}f_l^{kk}f_j^{ll}-f_j^{kk}f_k^{ll}f_l^{jj}) 
+f_j^{kk}f_l^{jj}(f_l^{ii}f_i^{kk}f_k^{ll}-f_i^{ll}f_k^{ii}f_l^{kk})\;.
\qquad}
$$
Under the $d$-semistability conditions 
the rank of the matrix is  5, and we obtain one infinitesimal
deformation, where the quadratic equations are perturbed. Furthermore
one has the perturbations of the cubic equations alone, which as before
have to be divisible by $x_i$. We already used $44$ 
coordinate transformations. The coefficient of $x_iy_\al y_\be$ can be
made to vanish with a transformation of the type $y_\ga \mapsto
y_\ga - \ep x_i$. So we have $28$ coefficients left and the diagonal
coordinate transformations, giving dimension 21.   
\qed

The computations in negative degree show that the tetrahedron $X$
is only a hyperplane section of threefolds with two-dimensional
singular locus, obtained by glueing together four cubic threefolds.

\neu %
We want to describe an explicit  deformation in the
$d$-semistable case.
We use the coordinate transformation $x_i\mapsto (f_4^{ii}/f_i^{44})x_i$,
$i=1$, $2$, $3$, 
which gives $f_i^{jj}f_4^{ii}/f_i^{44}$ as coefficient 
of $y_j^2$ in the new $F_i$. The $d$-semistability conditions yield that
the new coefficients satisfy
$f_i^{jj}=f_j^{ii}$. We will denote them by $f_{ij}$. A solution
to the linear equations above is then $d_{ij}=df_{kl}$, with
$d$ a new deformation variable.
Furthermore we use coordinate transformations to remove the $y_\al y_\be$
terms from the $F_i$.

We set $H_i=g_i^jy_j+g_i^ky_k+g_i^ly_l+h_ix_i$.
With this notation we get the following infinitesimal deformation:
$$
\displaylines{
x_ix_j+df_{kl}y_ky_l-df_{ij}f_{kl}(x_iy_j+x_jy_i)\;,\cr
x_iy_i+d(f_{kl}x_jH_j+f_{jl}x_kH_k+f_{jk}x_lH_l+
  f_{jk}f_{jl}y_j^2+f_{kl}f_{kj}y_k^2+df_{lj}f_{lk}y_l^2)\;, 
\vphantom{\bigoplus_a^a} \cr
\quad y_jy_ky_l-x_i(f_{ij}y_j^2+f_{ik}y_k^2+f_{il}y_l^2)-x_i^2H_i
\hfill\cr\hfill{}
  +dy_i((f_{ik}f_{jl}+f_{il}f_{jk})f_{ij}y_j^2+ 
  (f_{ij}f_{kl}+f_{il}f_{kj})f_{ik}y_k^2+ 
  (f_{ij}f_{lk}+f_{ik}f_{lj})f_{il}y_l^2)\; 
\hfill\cr\hfill{}
  +df_{ij}f_{ik}f_{il}y_i^3 
  +dy_i(f_{ik}f_{il}x_jH_j+f_{ij}f_{il}x_kH_k+f_{ij}f_{ik}x_lH_l) \;.\quad
}
$$

If we try to lift to higher order complicated formulas arise,
and it is not clear whether the computation is finite. It does stop 
if we restrict ourselves
to the case of tetrahedral symmetry. 
Then $f_{ij}$ does not depend
on $(i,j)$, and we call the common value $f$; likewise $g$ is the value
of all $g_i^j$ and $h$ of the $h_i$. We retain the notation
$H_i=g(y_j+y_k+y_l)+hx_i$. By a coordinate transformation
$x_i\mapsto x_i+df^2 y_i$ we simplify the expression for the first
6 equations. We write $t$ for the deformation parameter.

\proclaim Proposition.
The following set of equations defines a  degeneration of $K3$-surfaces
with special fibre a tetrahedron of degree 12:
$$
\displaylines{
x_ix_j+tf(y_k-tfH_k)(y_l-tfH_l)\;, \cr
x_iy_i+tf(x_jH_j+x_kH_k+x_lH_l)+tf^2(y_i^2+y_j^2+y_k^2+y_l^2)\;,  \cr
\qquad y_jy_ky_l-fx_i(y_i^2+y_j^2+y_k^2+y_l^2)-x_i^2H_i
 \hfill\cr\hfill{}
  -t^2f^2(y_jH_kH_l+y_kH_jH_l+y_lH_jH_k)+2t^3f^3H_jH_kH_l\;. \qquad
}
$$
The general fibre is a smooth $K3$-surface lying on
$\P^1\times\P^1\times\P^1$.

\roep Proof.
We suppress the computation. For $t\neq0$ the cubic equations
lie in the ideal of the quadrics. We can write three independent
equations 
$$
x_i(y_i-tfH_i)-x_j(y_j-tfH_j)
$$
which together with the first six define $\P^1\times\P^1\times\P^1$:
each square in the following picture gives an equation, where we put
$z_i:=\sqrt{-tf}(y_i-tfH_i)$.
$$
\setlength{\unitlength}{.6cm}
\begin{picture}(6,5.5)
\put(0.5,0){\line(1,0){3}}
\put(0,0.5){\line(0,1){3}}
\put(0.5,4){\line(1,0){3}}
\put(4,0.5){\line(0,1){3}}
\put(2.5,1){\line(1,0){1.4}}
\put(4.1,1){\line(1,0){1.4}}
\put(2,1.5){\line(0,1){2.4}}
\put(2,4.1){\line(0,1){0.4}}
\put(2.5,5){\line(1,0){3}}
\put(6,1.5){\line(0,1){3}}
\put(0.35,0.175){\line(2,1){1.3}}
\put(4.35,0.175){\line(2,1){1.3}}
\put(0.35,4.175){\line(2,1){1.3}}
\put(4.35,4.175){\line(2,1){1.3}}
\put(0,0){\makebox(0,0){$x_1$}}
\put(4,0){\makebox(0,0){$z_2$}}
\put(0,4){\makebox(0,0){$z_4$}}
\put(4,4){\makebox(0,0){$x_3$}}
\put(2,1){\makebox(0,0){$z_3$}}
\put(2,5){\makebox(0,0){$x_2$}}
\put(6,1){\makebox(0,0){$x_4$}}
\put(6,5){\makebox(0,0){$z_1$}}
\end{picture}
$$
We obtain the $K3$ by taking the complete intersection 
with the quadric
$\sum_i(x_iy_i+3tfx_iH_i+4tf^2 y_i^2)$.
\qed

\bfneu{The relation with quartics.}
We can construct a degree 12 tetrahedron from four planes in $\P^3$ 
by first blowing up each plane in 6 points and then glueing them
back together.

Therefore we first describe the blow up in a way adapted to our situation.
The problem of actually writing down a cubic equation is not often
treated in the literature (see however \cite{Ma}). 
Here we want to fix a particular hyperplane
section, which limits the choice of coordinates.   

Let two points lie on each side of the coordinate triangle in
$\P^2$ with coordinates $(z_1\cn z_2 \cn z_3)$. We describe them by
$z_k=z_i^2+a_{ij}z_iz_j+b_{ij}z_j^2=0$, where $(ijk)$ is a cyclic
permutation of $(123)$ (this means that we choose an orientation
on the triangle). 
As cubics through the six points we take
the coordinate triangle and three cubics, each consisting of a side
and a quadric passing through the remaining four points.
More precisely, we take
$$
\eqalign{
x_0&=z_1z_2z_3\cr
y_i&=z_i(z_k^2-a_{ki}z_kz_i+b_{ki}(z_i^2-a_{ij}z_iz_j+b_{ij}z_j^2))\;.
}
$$
One computes the relations
$$
z_k\>y_j-b_{ij}z_j\>y_k-
((1-b_{ij}b_{jk}b_{ki})z_i-a_{ij}z_j+b_{ij}b_{jk}a_{ki}z_k)\>x_0=0\;.
$$
By the Hilbert-Burch theorem the maximal minors of the relation
matrix 
$$
\pmatrix{
0&z_3&-b_{12}z_2&-(1-b_{12}b_{23}b_{31})z_1+a_{12}z_2-b_{12}b_{23}a_{31}z_3\cr
-b_{23}z_3&0&z_3&-(1-b_{12}b_{23}b_{31})z_2+a_{23}z_3-b_{23}b_{31}a_{12}z_1\cr
z_2&-b_{31}z_1&0&-(1-b_{12}b_{23}b_{31})z_3+a_{31}z_1-b_{31}b_{12}a_{23}z_2}
$$
give the cubics, up to a common factor $1-b_{12}b_{23}b_{31}$. 
By Lemma~(\ref{conditie}) this factor
vanishes exactly when the 6 points lie on a conic.

Viewing the relations as holding between the
$z_i$ gives the coefficient matrix
$$
\pmatrix{
(1-b_{12}b_{23}b_{31})x_0&b_{12}y_3-a_{12}x_0&b_{12}b_{23}a_{31}x_0-y_2 \cr   
b_{23}b_{31}a_{12}x_0-y_3&(1-b_{12}b_{23}b_{31})x_0&y_1b_{23}-a_{23}x_0  \cr  
b_{31}y_2-a_{31}x_0&b_{12}b_{31}a_{23}x_0-y_1&(1-b_{12}b_{23}b_{31})x_0  \cr
}\;.
$$
Its determinant is the equation of the surface.
After dividing by $b_{12}b_{23}b_{31}-1$ it equals
$$
\displaylines{
y_1y_2y_3-x_0(b_{12}y_3^2+b_{23}y_1^2+b_{31}y_2^2) 
\hfill\cr\hfill{}
-x_0^2(a_{12}a_{23}b_{31}y_2+a_{23}a_{31}b_{12}y_3+a_{31}a_{12}b_{23}y_1 
-(a_{12}y_3+a_{23}y_1+a_{31}y_2)(1+b_{12}b_{23}b_{31}))
\hfill\cr\hfill{}
-x_0^3((1-b_{12}b_{23}b_{31})^2+(1+b_{12}b_{23}b_{31})a_{12}a_{23}a_{31}
+b_{12}b_{23}a_{31}^2+b_{23}b_{31}a_{12}^2+b_{12}b_{31}a_{23}^2)\;.
}
$$
This last formula also works if the 6 points lie on a conic, but then it 
is easier to take the $y_i$ as product of a side and the conic through
the $6$ points; this means adding a multiple of $x_0$ to each $y_i$.
The equation then becomes $y_1y_2y_3-x_0Q(y)$ with $Q(z)$ the conic.

Now we apply this to our tetrahedron. We choose an orientation and 
orient the faces with the induced
orientation.
We get variables $x_i$ and $y_i$. For the face $i$
we take $x_i=z_jz_kz_l$ as before, but we multiply $y_j$ by a factor
$\la_{ij}$ to be determined later.  So we set $y_j=\la_{ij}z_j(z_l^2+\dots)$.
Now we look at the line $z_3=z_0=0$, with coordinates $(z_1\cn z_2)$.
Via the coordinates of  face $0$ we get the embedding
$(y_1\cn y_2)=(\la_{01} b_{31} z_1 \cn  \la_{02}  z_2 )$ whereas  face
$3$ gives $(y_1\cn y_2)=(\la_{31} z_1  \cn  \la_{32} b_{02} z_2 )$.
The condition that the Del Pezzo surfaces are glued in the same way
as the planes yields
the equations $\la_{01} \la_{32} b_{31} b_{02}=\la_{02}\la_{31} $.
By even permutations of $(0123)$ we get in total six equations.
They are solvable if and only if
$$
b_{01}b_{10}b_{02}b_{20}b_{03}b_{30}b_{12}b_{21}b_{13}b_{31}b_{23}b_{32}=1\;,
$$
a condition obtained by multiplying the six equations.

The $d$-semistability conditions 
$
f_k^{jj}f_l^{kk}f_j^{ll}-f_j^{kk}f_k^{ll}f_l^{jj}=0 
$
give for the cubic above, using $(jkl)=(123)$:
$$
{b_{30}\over \la_{21}^2}
{b_{10}\over \la_{32}^2}
{b_{20}\over \la_{13}^2}
=
{b_{03}\over \la_{12}^2}
{b_{01}\over \la_{23}^2}
{b_{02}\over \la_{31}^2}\;.
$$
Using that the $\la_{ij}$ satisfy the equations 
$\la_{32}b_{02}/ \la_{31}=\la_{02}/(\la_{01} b_{31})$
we see that this condition is equivalent to
$b_{21}^2b_{32}^2b_{13}^2=b_{12}^2b_{23}^2b_{31}^2$, which is one of
the conditions that the 24 points are cut out by a quartic.

We can ask  which choices of 24 points
give our symmetric tetrahedron. The condition $\prod b_{ij}=1$
limits the possibilities. In particular, if all $b_{ij}=1$, the 
six points in each face lie on a conic, giving a singular tetrahedron. 
If we take the quartic
$Q=(a \sigma_1^2  + b \sigma_2)^2$ then each element of the pencil
has 12 singular points. We can blow up them and blow down the six 
conics in the faces by embedding the pencil in $\P^7\times \P^1$
with the linear system of cubics in $\P^3$ with as base points the 12 singular
points. We set 
$$
\eqalign{
x_i&=z_jz_kz_l\cr
y_i&=z_i(a \sigma_1^2  + b \sigma_2)\;.
}
$$
We obtain a symmetric tetrahedron with $g=h=0$.
 
We get nonsingular Del Pezzo surfaces by taking all $b_{ij}=-1$,
and $a_{ij}=a$. Then  $f=-1$, $g=-a^2$ and $h=a^2+4$.
The points on the side of the tetrahedron are given by
$$
(z_i^2+az_iz_j-z_j^2)(-z_i^2+az_iz_j+z_j^2)=(-z_i^4+(2+a^2)z_i^2z_j^2-z_j^4)\;.
$$

In particular, we obtain different smoothings of the same tetrahedron,
those embedded in $\P^7$ and other ones, where the general fibre
is embeddable in $\P^3$. They belong to different 19-dimensional
hypersurfaces in the 20-dimensional subspace of the versal
deformation whose general fibre is a smooth $K3$-surface.

\section{Deformation theory}
%
\neu
Let $X=\cup X_i$ be a normal crossings surface with normalisation
$\wtx = \coprod  X_i$.
The components of the double locus $D$ are $D_{ij}=X_i\cap X_j$.
The divisor $D_i:=\cup_j D_{ij}$ is a normal crossings divisor in $X_i$.
We set $\wl D= \coprod D_i$.

As $X$ is locally a hypersurface in a 3-fold $M$, 
its cotangent cohomology sheaves $\ct_X^i$
vanish for $i\geq2$ and 
$$
0\lra \ct_X^0 \lra \Theta_M|_X \lra N_{X/M} \lra \ct_X^1 \lra 0\;.
$$
There is a canonical isomorphism 
$ \ct_X^1 \cong \sier D(X)$ and in particular, if $X$ is $d$-semistable,
then $ \ct_X^1 \cong \sier D$ \cite[Prop.~2.3]{Fr}.

\proclaim Lemma.
There is an exact sequence
$$
0\lra \ct^0_X \lra n_*\Theta_{\wtx} (\log \wl D)\lra  \ct^0_D \lra 0\;.
$$

\roep Proof.
This is a local computation. The  sheaf $\Theta_M(\log X)$ of 
vectors fields on $M$ which preserve
$z_1z_2z_3=0$ is generated by the $z_i{\partial \over \partial z_i}$.
Restricted to a component $X_i\colon z_i=0$ we get sections of
$\Theta_{X_i}(\log D_i)$. The restrictions to different components
satisfy the obvious compatibility condition.
\qed

Sections of $\ct^0_D$ are given by vectorfields on each component,
which vanish in the triple points. We study $\Theta_{X_i}(\log D_i)$
with the exact sequence
$$
0 \lra \Theta_{X_i}(\log D_i) \lra \Theta_{X_i} \lra \oplus_j
N_{D_{ij}/X_i} \lra 0 \;.
$$
For a $d$-semistable $K3$-surface $X$ in $(-1)$-form  
$$H^0(D_{ij},N_{D_{ij}/X_i}) =
H^1(D_{ij},N_{D_{ij}/X_i}) = 0\;.$$
Each component  
$X_i$ is $\P^2$ blown up in $k\geq 3$ points and $H^2(\Theta_{X_i})=0$,
$h^0(\Theta_{X_i})=\max(0,8-2k)$, $h^1(\Theta_{X_i})=\max(0,2k-8)$.

So  $H^0(\Theta_{X_i})\neq 0$ only in the case that $k=3$ and 
the double curve $D_i$ is
a hexagon. We then call $X_i$ a hexagonal component, or hexagon for short.

\proclaim Lemma {\rm \cite[Cor.~3.5]{Frr}}. \label{verdwijn}%
For a $d$-semistable $K3$-surface $X$ of type III
in $(-1)$-form   $H^0(X, \ct_X^0)=0$.

\roep Proof.
We first describe the sections of $H^0(\Theta_{X_i})$
for a hexagonal component. We blow up $\P^2$ in the vertices
of the coordinate triangle.
As basis for the linear system of cubics we take the monomials
given by black dots in the picture below.

$$
\setlength{\unitlength}{1cm}
\begin{picture}(3.8,3.8)(-.4,-.4)
\put(0,0){\circle{.2}}
\put(-.2,-.2){\makebox(0,0)[tr]{$z_1^3$}}
\put(1,0){\circle*{.2}}
\put(1,-.2){\makebox(0,0)[t]{$x_6$}}
\put(2,0){\circle*{.2}}
\put(2,-.2){\makebox(0,0)[t]{$x_5$}}
\put(3,0){\circle{.2}}
\put(3.2,-.1){\makebox(0,0)[tl]{$z_3^3$}}
\put(.5,.86){\circle*{.2}}
\put(.3,.86){\makebox(0,0)[r]{$x_1$}}
\put(1.5,.86){\circle*{.2}}
\put(1.5,.66){\makebox(0,0)[t]{$x_0$}}
\put(2.5,.86){\circle*{.2}}
\put(2.7,.86){\makebox(0,0)[l]{$x_4$}}
\put(1,1.73){\circle*{.2}}
\put(.8,1.73){\makebox(0,0)[r]{$x_2$}}
\put(2,1.73){\circle*{.2}}
\put(2.2,1.73){\makebox(0,0)[l]{$x_3$}}
\put(1.5,2.6){\circle{.2}}
\put(1.5,2.8){\makebox(0,0)[b]{$z_2^3$}}
\end{picture}
$$

A vector field $\th$ on $X_i$ comes from a vector field on $\P^2$
which vanishes in the points blown up.
We can give it homogeneously by $a_1z_1\pab{z_1}+a_2z_2\pab{z_2}
+a_3z_3\pab{z_3}$, subject to the relation $z_1\pab{z_1}+z_2\pab{z_2}
+z_3\pab{z_3}=0$. In the $x_j$ coordinates we get
$$\displaylines{\textstyle\quad
(a_1+a_2+a_3)x_0\pab{x_0}+(2a_1+a_2)x_1\pab{x_1}+
(2a_2+a_1)x_2\pab{x_2}+(2a_2+a_3)x_3\pab{x_3}
\hfill\cr\hfill\textstyle{}+(2a_3+a_2)x_4\pab{x_4}+
(2a_3+a_1)x_5\pab{x_5}+(2a_1+a_3)x_6\pab{x_6}\;.\quad}
$$
We restrict to the line $(x_{j-1}\cn x_{j})$ and take as generator
of $\ct^0_D|_{D_{ij}}$ the vector field $\th_j=\half (x_j\pab{x_j}-
x_{j-1}\pab{x_{j-1}})$. On the $(x_6\cn x_1)$-line 
$\th = (a_2-a_3) \th_1$ and on the $(x_1\cn x_2)$-line
$\th = (a_2-a_1) \th_2$. The remaining coeffcients $\th = \be_j \th_j$
are found by cyclic permutation. They satisfy $\be_j=\be_{j-1}+
\be_{j+1}$. In particular, two adjacent coefficients determine all 
other ones and opposite coefficients add up to zero.

\begin{figure}
\centering
\setlength{\unitlength}{1.5cm}
\begin{picture}(3,4)(-1,0)
\put(0,0){\line(1,1){1}}
\put(0,2){\line(1,-1){1}}
\put(0,2){\line(1,1){1}}
\put(0,4){\line(1,-1){1}}
\put(0,2){\line(-1,0){1}}
\put(1,1){\line(1,0){1}}
\put(1,3){\line(1,0){1}}
\put(1.5,0){\makebox(0,0)[lb]{$X_0$}}
\put(1.5,2){\makebox(0,0)[l]{$X_1$}}
\put(-.5,1){\makebox(0,0)[r]{$X_1'$}}
\put(1.5,4){\makebox(0,0)[lt]{$X_2$}}
\put(-.5,3){\makebox(0,0)[r]{$X_2'$}}
\put(1.5,1.1){\makebox(0,0)[b]{$0$}}
\put(.4,.6){\makebox(0,0)[rb]{$0$}}
\put(1.5,2.9){\makebox(0,0)[t]{$0$}}
\put(1.5,3.1){\makebox(0,0)[b]{$0$}}
\put(.6,1.6){\makebox(0,0)[lb]{$\be$}}
\put(.6,2.4){\makebox(0,0)[lt]{$\be$}}
\put(-.5,2.1){\makebox(0,0)[b]{$\be$}}
\put(.4,1.4){\makebox(0,0)[rt]{$-\be$}}
\put(.4,2.6){\makebox(0,0)[rb]{$-\be$}}
\put(-.5,1.9){\makebox(0,0)[t]{$-\be$}}
\put(.4,3.4){\makebox(0,0)[rt]{$-2\be$}}
\put(.6,3.6){\makebox(0,0)[lb]{$2\be$}}
\end{picture}
\caption{Coefficients of a vector field.}
\label{coefs}
\end{figure}
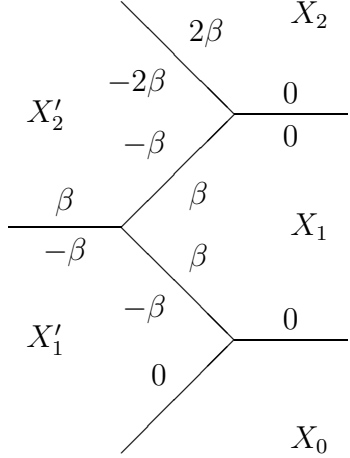

Let $\th\in H^0(X, \ct_X^0)$ be
a non-vanishing global section.
As the dual graph is a triangulation of $S^2$ one has
$\sum_i (6-e_i)=12$, where $e_i$ is the number of components
of the double curve $D_i$. So there exist non-hexagonal components, and
$\th$ vanishes on them.
Suppose $\th$ vanishes on $X_0$ and not on the adjacent hexagon $X_1$.
We are going to look at the restriction of $\th$ to other components,
as illustrated in Fig.~\ref{coefs}.

Let $T=X_0\cap X_1 \cap X_1'$ be a triple point. We know
that $\th$ vanishes on $X_1\cap X_0$. If it also vanishes
on $X_1\cap X_1'$, then it vanishes
altogether, contrary to the assumption. Therefore $X_1'$ is also hexagonal.
Let $\th= \be \th_0$ on $D_{ii'}=X_1\cap X_1'\subset X_1$.
Considered on $X_1'$ the restriction of $\th$ is $-\be$ times the
generator. 
The other triple point on $D_{ii'}$ involves a hexagon $X_2'$, which 
contains also the triple point $X_1\cap X_2' \cap X_2$. Considered
on $X_2'$, the coefficient of the restriction of $\th$ 
to $X_2'\cap X_1'$ is $\be$, to $X_2'\cap X_1$ it is $-\be$,
so to $X_2'\cap X_2$ it is $-2\be$. Therefore on $X_2$ $\th$ has 
adjacent  coefficients $0$, $2\be$.
Inductively we find components $X_n'$, $X_n$ with the coefficient
$n\be$ occurring. As there are only finitely many components, this is 
impossible.
\qed


\proclaim Theorem. \label{tangent}
Let $X=\cup_{i=1}^k X_i$ be
a $d$-semistable $K3$-surface  of type III
in $(-1)$-form, with $k$ components. 
Then 
$$
\eqalign{
\dim H^1(X, \ct_X^0)&= k+18 \cr
\dim H^0(X, \ct_X^1)&= 1 \cr
\dim H^1(X, \ct_X^1)&= k-1 }
$$
So $\dim T^1_X=k+19$, $\dim T^2_X=k-1$.

\roep Proof.
As the dual graph triangulates $S^2$ we have $V-E+F=2$, where $V=k$,
the number of components of $X$, $E$ is the number of double curves
and $F$ is the number of triple points. Each double curve contains
two triple points, so  $F=\break 2/3 E$, which makes $E=3k-6$.
A component $X_i$, which is $\P^2$ blown up in $\delta_i$ points, has 
$e_i=9-\delta_i$ double curves. Observe that $\sum_i e_i= 2E$.
The exact sequence above gives $\dim H^1(X, \ct_X^0)= \sum_i 2(5-e_i)+E=
10 V - 3 E=k +18$.

We have $h^0(X, \ct_X^1)=h^0(D,\sier D)=1$ and $h^1(X, \ct_X^1)=h^1(D,\sier D)=
1-\chi=1-(E-2F)=k-1$.
\qed

\neu\label{obstruc}
Locally trivial deformations of a $d$-semistable $K3$-surface $X$
are unobstructed and fill up a codimension one smooth subspace
of the base of the versal deformation with tangent space
$H^1(X,\ct^0_X)$. This means that every equation of the base is
divisible by the equation of this hypersurface.
As one obtains the base space as fibre of a map 
$T^1 \to T^2$, we look at the 
map 
$$
{\rm Ob}\colon H^1(\ct_X^0)\times H^0(\ct_X^1) \to H^1(\ct_X^1)\;.
$$
Let $\xi$ be a global generator of  $\ct_X^1$. 
The existence of a second smooth component (of dimension 20) follows,
if one can show that 
${\rm Ob}(\;.\;,\xi)\colon H^1(\ct_X^0) \to H^1(\ct_X^1)$ 
is a surjective linear map.
To describe it we start with the 
map ${\rm Ob}(\;.\;,\xi)\colon \ct_X^0 \to \ct_X^1$.
Locally $X$ is a hypersurface given by an equation $f=0$ and elements
of $\ct^0_X$ come from ambient vector fields satisfying $\th(f)=c\,f$.
We can choose coordinates such that $\xi$ acts as $f \mapsto 1$.
Then ${\rm Ob}(\th,\xi)=-c\,\xi$. In the normal crossings situation
the map ${\rm Ob}(\;.\;,\xi)$ is surjective
and we get an exact sequence
$$
0 \lra {\cal S} \lra \ct_X^0 \lra \ct_X^1 \lra 0 \;.
$$

The kernel of the map
${\rm Ob}(\;.\;,\xi)\colon H^1(\ct_X^0) \to H^1(\ct_X^1)$
can be characterised in a different way (\cite{FS}).
If $X=\cup_{i=1}^k X_i$  occurs as central fibre in 
a degeneration ${\cal X} \to S$, we  define $k$ line bundles
$L_i:=\sier{\cal X}(X_i)|_X$. On a $d$-semistable $X$ it can be defined by
$$
\eqalign{
L_i|_{X_i}&=\sier {X_i}(-D_i) \cr
L_i|_{X_j}&=\sier {X_j}(X_i\cap X_j), \qquad j\neq i
}
$$
with appropriate glueings, using the global section of $\sier D(X)$.
The bundle $L_i$ defines a class $\xi_i$ in
$$
H^1(X,\sier X^*)\cong \ker\left\{ H^2(X,\Z)\to H^2(\sier X)=\C\right\}\;,
$$
which therefore lies in $H^1(\Omega^1/\tau^1)$, where $\Omega^1/\tau^1$ 
are the K\"ahler differentials modulo torsion \cite[Sect.~1]{Fr}.
The condition that $L_i$ lifts to line bundles on a
locally trivial deformation with tangent vector $\th \in H^1(\ct_X^0)$
is that $\langle\th,\xi_i\rangle=0$ with $\langle-,-\rangle$
the perfect pairing 
$H^1(\ct_X^0) \otimes H^1(\Omega^1/\tau^1) \to H^2(\sier X)=\C$ 
\cite[(2.10)]{Fr}.
The surjectivity of the map
${\rm Ob}(\;.\;,\xi)$ follows from the following lemma.

\proclaim Lemma. \label{span}%
The classes $\xi_i$ span a $(k-1)$-dimensional
subspace of $H^2(X,\Z)$.

\roep Proof.
We compute $H^2(X,\Z)$ as the kernel of the map
$$
\oplus H^2(X_i,\Z) \to \oplus H^2(D_{ij},\Z)\;.
$$
Each $\xi_i$ gives rise to a divisor $\sum_m a_{lm}D_{lm}$ on $X_l$,
$l=1$, \dots, $k$, with coefficients satisfying $a_{lm}+a_{ml}=0$
(and $a_{lm}\neq0$ only if $i=l$ or $i=m$). The relation
$\sum \xi_i=0$ holds.

Let now $\sum b_i\xi_i =0\in H^2(X,\Z)$. It gives rise to a divisor 
$\sum_m \be_{lm}D_{lm}$ on $X_l$. If the classes $D_{lm}$ are independent
in $H^2(X_l,\Z)$, then  $\be_{lm}=0$ for all $m$. 
This condition is not satisfied if $X_l$ is a hexagon.
Then we can only conclude that $\be_{l,m-1}+\be_{l,m+1}=\be_{lm}$.
With the same argument as in the proof of Thm.~(\ref{verdwijn}),
illustrated by Fig.~\ref{coefs}, we infer that even in this case
$\be_{lm}=0$ for all $m$. 

Therefore $b_i=b_j$ for all pairs $(i,j)$ such that $X_i\cap X_j
\neq \emptyset$. This implies that $\sum b_i\xi_i$ is a multiple of 
$\sum \xi_i$.
\qed

We summarise:

\proclaim Theorem {\rm \cite[(5.10)]{Fr}}.
A $d$-semistable $K3$-surface $X$ of type III is smoothable.
Its versal base space is the union $V_1\cup V_2$, where
$V_1$ is a smooth hypersurface corresponding
to locally trivial deformations of $X$, which meets transversally
a $20$-dimensional smooth subspace $V_2$, with $V_2\setminus V_1$
parametrising smooth $K3$-surfaces and $V_2\cap V_1$ 
locally trivial deformations of $X$ for which $\sier D(X)$ remains trivial.

\bfneu{Embedded deformations.} 
We relate the above  results to direct computations with generators
and relations for the cone over $X$, as for the tetrahedron.
The case of cones  over non-singular varieties is treated in \cite{St}. 
We suppose that the affine cone $C(X)$ over $X$ is Cohen-Macaulay.
The starting point is the exact sequence
\begin{equation}
0\lra T^0_{C(X)} \lra \Theta_{\C^{n+1}}|_{C(X)} \lra N_{C(X)} 
\lra T^1_{C(X)} \lra 0\;,\label{rij}
\end{equation}
which we shall relate to exact sequences of sheaves on $X$.
We set $U=C(X)\setminus0$; then $\pi\colon U \to X$  
is a $\C^*$-bundle over $X$.
For a reflexive sheaf $\cf$ on $C(X)$ we have
$H^0(C(X),\cf)=H^0(U,\cf)$. All sheafs $\cf$ considered here have a 
natural $\C^*$-action, so
$\pi_*\cf$ decomposes into the direct sum of eigenspaces. In particular,
the degree 0 part is the sheaf of $\C^*$-invariants.
With homogeneous coordinates $x_i$ the $\C^*$-invariant
sections $x_j\pab{x_i}$ of $H^0(U,\Theta_{\C^{n+1}}|_{C(X)})$ can be 
considered as elements of $H^0(X, V^*\otimes_\C \sier X(1))$,
where $V=H^0(X,\sier X(1))$.
We get the degree zero part $T^1_{C(X)}(0)$ as 
$\coker  H^0(X, V^*\otimes_\C \sier X(1)) \to H^0(X, N_{X/\P^n})$.
We factorize this map corresponding  to a splitting of the
exact sequence (\ref{rij}):
\begin{eqnarray}
0\lra T^0_{C(X)} \lra \Theta_{\C^{n+1}}|_{C(X)} \lra G  \lra 0
\label{eerste}\\
0\lra G \lra N_{C(X)} \lra T^1_{C(X)} \lra 0\;.\label{tweede}
\end{eqnarray}
Denoting by $\cg_X$  the sheaf of $\C^*$ invariants associated to $G$
we obtain 
$$
 H^0(X, V^*\otimes_\C \sier X(1)) \lra H^0(X,\cg_X) \lra H^0(X, N_{X/\P^n})\;.
$$
On $X$ we have the exact sequence
$$
0\lra \cg_X \lra N_{X/\P^n} \lra \ct^1_{X} \lra 0\;.
$$
The short exact sequence (\ref{eerste}) gives
$$
0\lra \D_X \lra V^*\otimes_\C \sier X(1) \lra \cg_{X} \lra 0
$$
with $\D_X$ the sheaf of differential operators on $X$, which is related
to $\ct^0_X$ by the exact sequence
$$
0 \lra \sier X \lra \D_X \lra \ct^0_X \lra 0\;.
$$

\proclaim Proposition.  \label{ingebed}%
Let $X$ be  a $d$-semistable $K3$ of type III in $(-1)$-form.
The space of infinitesimal locally trivial embedded deformations 
is $H^1(X,\D_X)$, of dimension $k+17$. It has codimension one in 
$T^1_{C(X)}(0)$.

\roep Proof.
From the computation of $h^i(\ct^0_X)$ in 
(\ref{tangent}) and the exact sequence for $\D_X$
we conclude that $h^0(\D_X)=h^0(\sier X)=1$. As $h^1(\sier X)=0$
and $h^2(\sier X)=1$ we get the exact sequence
$$
0 \lra  H^1(X,\D_X) \lra H^1(X,\ct^0_X) \lra H^2(X,\sier X)
\lra H^2(X,\D_X) \lra 0 \;.
$$
The  line bundle $\sier{}(1)$ determines a
class $h\in H^1(\Omega^1/\tau^1)$, which lifts to a deformation 
$\th\in H^1(X,\ct^0_X)$ if and only if 
$\langle\th,h\rangle=0$ with $\langle-,-\rangle$
the perfect pairing 
$H^1(\ct_X^0) \otimes H^1(\Omega^1/\tau^1) \to H^2(\sier X)=\C$.
This accounts for the non-algebraic deformation direction.
So $\dim H^1(\D_X)=k+17$ and $H^2(\D_X)=0$.
We then obtain 
$$
H^1(X,\D_X)= \coker \left\{H^0(X, V^*\otimes_\C \sier X(1)) 
\lra H^0(X,\cg_X)\right\}
$$ 
and $h^1(\cg_X)=0$, as $h^i(X,\sier X(1))=0$ for $i>0$.
Finally we get $H^1(N_{X/\P^n})=H^1(\ct^1_{X})$ and
the exact sequence
$$
0 \lra H^0(X,\cg_X) \lra H^0(X,N_{X/\P^n}) \lra H^0(X,\ct^1_X)
\lra 0\;.
$$
\qed

\neu
For $T^2_{C(X)}(0)$ we can argue as in the smooth case
\cite[(1.25)]{St} to obtain the exact sequence
$$
0 \lra T^2_{C(X)}(0) \lra
H^1(X,N_{X/\P^n}) \lra \bigoplus H^1(X,\sier X(d_j))
$$
with the $d_j$ the degrees of the generators of the ideal of $C(X)$ (or
of $X$). In particular, in our situation
$ T^2_{C(X)}(0)= H^1(N_{X/\P^n})=H^1(\ct^1_{X})$.

\proclaim Theorem {\rm \cite[(5.5)]{FS}}.
A $d$-semistable $K3$-surface $X$ of type III in $\P^n$ is smoothable
by embedded deformations.
They form a $19$-dimensional smooth component.

\roep Proof.
In the embedded case the base space is also the fibre
of a map between the relevant cotangent modules, and the 
locally trivial deformations are unobstructed.
The map ${\rm Ob}\colon H^1(\D_X)\times H^0(\ct_X^1) \to H^1(\ct_X^1)$ is the
restriction of the obstruction map in (\ref{obstruc}). We observe  that
$H^1(\D_X)$ is transversal to 
$\cap_i\ker {\rm Ob}(\;.\;,\xi_i)$, as the class $h$ satisfies $h^2>0$
and is therefore independent of the classes of the $\xi_i$.
\qed

\bfneu{The topology of the special fibre.}
One can compute the homology $H_*(X,\Z)$
with a Mayer-Vietoris spectral sequence
\cite[Prop.~2.5.1]{Pe} with $E^1$-term $E^1_{p,q}=H_p(X^{[q]},\Z)$,
where $X^{[0]}=\coprod X_i$, $X^{[1]}=\coprod D_{ij}$ and $X^{[2]}$
the set of triple points $P_{ijk}=X_i\cap X_j \cap X_k$.

\proclaim Proposition.
Let $X=\cup_{i=1}^k X_i$ be
a $d$-semistable $K3$-surface  of type III
in $(-1)$-form, with $k$ components. 
Then 
$$
\eqalign{
\dim H^0(X, \Z)&= 1 \cr
\dim H^2(X, \Z)&= k+19 \cr
\dim H^4(X, \Z)&= k }
$$

\roep Proof.
The $E^1$-term of the spectral sequence looks like
$$
\def\normalbaselines{\baselineskip18pt
 \lineskip2pt \lineskiplimit2pt}
\matrix{
\oplus H_4(X_i,\Z) \cr
0\cr
\oplus H_2(X_i,\Z) & \oplus H_2(D_{ij},\Z) \cr
0&0\cr
\oplus H_0(X_i,\Z) & \oplus H_0(D_{ij},\Z) & \oplus H_0(T_{ijk},\Z) }
$$
To prove that the map 
$\oplus H_2(D_{ij},\Z) \to \oplus H_2(X_i,\Z)$ is injective
we observe that $\oplus_j H_2(D_{ij},\Z) \to H_2(X_i,\Z)$ is
injective unless $X_i$ is a hexagonal component.
We take care of those by argueing as in the proofs
of Lemmas (\ref{verdwijn}) and (\ref{span}). If the component
$X_i$ is obtained by blowing up $\P^2$ in $\delta_i$ points,
then $b_2(X_i)=\delta_i+1=10-e_i$ with the notation of
(\ref{verdwijn}), so the cokernel of the map 
$\oplus H_2(D_{ij},\Z) \to \oplus H_2(X_i,\Z)$
has dimension $10V-3E=k+18$. The dimension formulas now follow
from the spectral sequence.
\qed

\neu \label{klasse}%
We describe the non-algebraic homology class in more detail.
Each double curve contains two triple points, which are homologous,
so the boundary of an interval. On a component $X_i$ these intervals 
make up a closed polygon (with $e_i$ edges), which itself is
the boundary of a topological disc.
For the case of $\P^2$ blown up in 4 points this is illustrated in 
Fig.~\ref{blaas}: after blowing up we have a pentagon, which is the boundary
of the strict transform of the shaded area. With the given coordinates
this strict transform consists of all points on the Del Pezzo surface
with positive coordinates.
Finally the discs glue together to a real polyhedron with the same
dual graph as the complex surface $X$. 

\neu
A nice construction for studying the homology
the general fibre  is given by \cite{ACa}.
Let $\sigma_i\colon {\cal Z}_i \to \cal X$ be the oriented real
blow-up of $X_i\subset \cal X$. This is a manifold with
boundary, whose boundary $\partial {\cal Z}_i= \sigma_i^{-1}(X_i)$
is isomorphic to the boundary of a tubular neighbourhood of $X_i$ in
$\cal X$. The fibred product $\sigma\colon {\cal Z} \to \cal X$
of the $\sigma_i$ is a manifold with corners. Its boundary ${\cal N}
:=\partial {\cal Z}$ comes with a map to $X$. It also fibres over
$S^1$: the composed map ${\cal Z} \to {\cal X} \to S \ni 0$ 
extends to a map from  ${\cal Z}$ to the real oriented blow-up
of $S$ in 0 (polar coordinates!).
A  fibre of ${\cal N} \to S^1$ is then a topological model of the
general fibre.

This model is not sufficient to describe the monodromy. One has first to
replace $X$ by the geometric realisation of the simplicial
object $X^{[\cdot]}$: one replaces each double point by an interval,
and each triple point by a 2-simplex. A final fibred product then gives
the new model. For details see \cite[\S 2]{ACa}.
 
\section{Hodge algebras}

\bfneu{Stanley--Reisner rings.}
%
%
Let $\De$ be a simplicial complex with set of vertices 
$V=\{v_1, \dots, v_n\}$. A monomial on $V$ is an element of
$\N^V$. Each subset of $V$ determines a
monomial on $V$ by its characteristic function.
The support of a monomial $M\colon V \to \N$ is the set
$\supp M=\{v\in V \mid M(v)\neq 0\}$.
The set $\Sigma_\De$ of monomials whose support is not a face
is an ideal, generated by the monomials corresponding to minimal 
non-simplices.

Given a ring $R$ and an injection $\phi\colon V \to R$
we can associate to each monomial $M$ on $V$ the element
$\phi(M)=\prod _{v\in V} \phi(v)^{M(v)}\in R$. We will usually identify
$V$ and $\phi(V)$ and write $M\in R$ for $\phi(M)$. This applies in
particular to the polynomial ring $K[V]$ over a field $K$. 
The ideal $\Sigma_\De$
gives rise to the  Stanley-Reisner ideal
$I_\De\subset K[V]$. The {\sl Stanley-Reisner ring\/} is
$A_\De= K[V]/I_\De$.

Deformations of Stanley-Reisner rings are studied in \cite{AC}.
  
\nroep Example.
Let $\De$ be an octahedron. We map the set of vertices to $\C[x_1,\dots,x_6]$
such that
opposite vertices correspond to variables with index sum 7.
\begin{figure}[h]
\centering
\includegraphics[width=6cm]{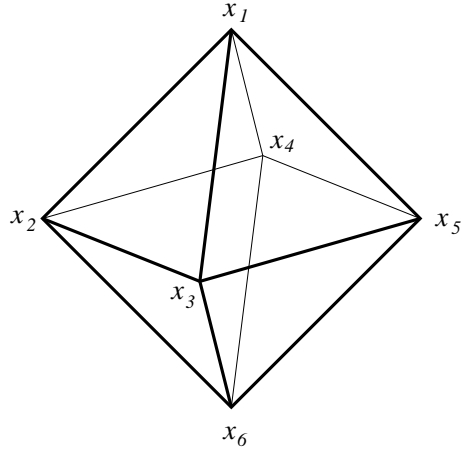}
\caption{Octahedron.}
\end{figure}
The Stanley-Reisner ring is minimally generated by  the three monomials
$x_ix_{7-i}$. The spaces smoothes to a $K3$-surface, the complete
intersection of three general quadrics.
A general 1-parameter deformation is not semi-stable, because the
total space has singularities at the six quadruple points of the special
fibre. 
\endroep

\nroep Remark. To get an  octahedron as dual graph we need the 
incidence relations of a cube.
The toric variety associated to a cube is $\P^1\times
\P^1\times \P^1$. The general anticanonical divisor is a smooth 
$K3$, whereas the complement of the torus is the unions of six
quadrics (each of the form $\P^1\times \P^1$) with dual graph the octahedron.
A small resolution of a general pencil yields a
semistable degeneration.

\nroep Definition {\rm \cite{DEP}}.
Let $H$ be a finite partially ordered set and $\Sigma$
an ideal of monomials on $H$. Let $K$ be a commutative ring and $R$ 
a commutative $K$-algebra. 
Suppose  an injection $\phi\colon H \to R$ is given.
Then $R$ is a {\sl Hodge algebra\/} (or  algebra with straightening
law) {\sl governed by $\Sigma$\/} if
\begin{itemize}
\item[H-1] $R$ is a free $K$-module admitting the set of
monomials not in $\Sigma$ as basis
\item[H-2] for each generator $M$ of $\Sigma$ in the unique expression
$$
M=\sum_{N\notin \Sigma} k_{M,N} N, \qquad k_{M,N}\in K\;, \eqno(*)
$$  
guaranteed by H-1, for each $x\in H$ dividing $M$ and each
$N\notin \Sigma$ with $k_{M,N}\neq0$  there is a $y_{M,N}$ dividing
$N$ and satisfying $y_{M,N}<x$.                                                
\end{itemize}

\noindent 
The relations $(*)$ are called {\sl straightening relations\/} for $R$.

\neu 
If $R$ is graded and the elements of $\phi(H)$ are homogeneous
the straightening relations give a presentation for $R$
\cite[p.~15]{DEP}.

We note that $R$ is a deformation of the discrete Hodge algebra
governed by $\Sigma$, whose ideal is generated by the monomials $M$.

\nroep Example. The equations of the tetrahedron of degree 12 
of (\ref{tetra}) are straightening relations. We take $\Sigma$ as
Stanley-Reisner
ideal $\Sigma_\De$, where $\De$ is the stellation of the tetrahedron:
in each top-dimensional face we take an additional vertex, which is joined
to all vertices on the face. The partial order on the set of vertices  
is obtained by declaring the new vertices to be smaller.
The discrete Hoge algebra has then equations
$x_ix_j$, $x_iy_j$ and $y_jy_ky_l$.

\section{The dodecahedron}

\neu
To get an icosahedron as dual graph we need the incidence relations
of a dodecahedron. Each side should be a rational surface and the intersection
with the other surfaces should have a pentagon as dual graph.
A pentagon occurs as hyperplane section of a Del Pezzo surface of degree
5. So we can realise our dodecahedron by glueing together 12 
Del Pezzo surfaces.

\begin{figure}[h]
\centering
\begin{minipage}[c]{8cm}
  \includegraphics[width=\textwidth]{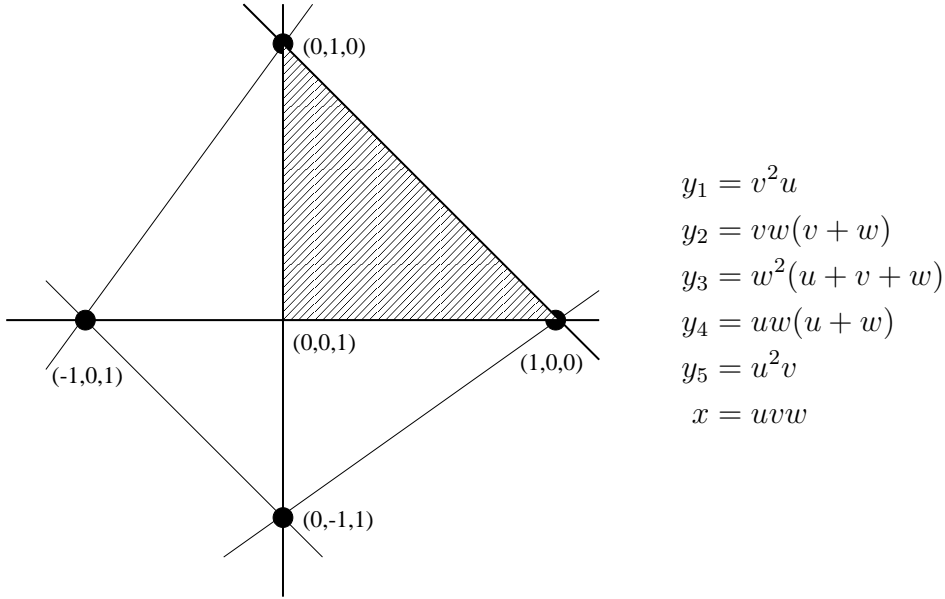}
\end{minipage}
$\qquad\eqalign{ 
y_1&=v^2u \cr
y_2&=vw(v+w) \cr
y_3&=w^2(u+v+w) \cr
y_4&=uw(u+w) \cr
y_5&=u^2v\cr
x&=uvw 
}$
\caption{Blowing up $\P^2$ in 4 points.}\label{blaas}
\end{figure}

We first describe the Del Pezzo surfaces. Each of those is an
extension of its pentagonal hyperplane section. 
Its coordinate ring can be obtained as
Stanley-Reisner ring of a pentagon as 1-dimensional simplicial complex.
Introducing variables $y_i$, we get the equations $y_{i-1}y_{i+1}$.
With an extra variable $x$ the Del Pezzo surface 
has equations 
$$
y_{i-1}y_{i+1}-xy_i-x^2\;.
$$
These are the pfaffians of the matrix
$$
\pmatrix{
 0  &  y_1  &  x  &  -x  & -y_5 \cr
- y_1 & 0  & y_2 &  x  &  -x    \cr
-x & - y_2 & 0  & y_3 &  x      \cr
x  &  -x & - y_3 & 0  & y_4      \cr
y_5 & x  &  -x & - y_4 & 0  }\;.
$$
We can check that this is indeed a smooth Del Pezzo of degree 5
by giving an explicit birational map from $\P^2$, which blows up
four points, see Fig.~\ref{blaas}.
To the variable $x$ corresponds a new vertex  at the centre of the
pentagon. By joining it to all other vertices  we obtain 
a  2-dimensional simplicial complex, 
and the homogeneous coordinate ring of the Del Pezzo surface
is a graded Hodge algebra  governed by the Stanley-Reisner ideal
of the complex:
to satisfy H-2  we take $x$ to be less than all $y_i$.

\neu
To construct a normal crossings dodecahedron of degree 60
we glue twelve Del Pezzo surfaces. 
We get a simplicial complex $\Dee$ by stellating
a dodecahedron: we take in each face the center of the pentagon as extra
vertex. A non-convex realisation of this complex is the
great stellated dodecahedron. 

\proclaim Proposition.
The coordinate ring of the dodecahedron of degree 60
is a  graded Hodge algebra governed by $\Sigma_\Dee$.

We describe the equations for a dodecahedron $X$ 
with icosahedral symmetry in more detail.
We have   
$20$ variables $y_\al$, one for each dodecahedral vertex 
and $12$ variables $x_i$ from the extra vertices  in the faces. 
We will denote the vertices  
by $\al$ and $i$.
As two face vertices are not connected by an edge we have
66 equations $x_ix_j=0$, $i\neq j$. If $\overline{i\al}$
is not an edge, we have $x_iy_\al=0$; there are $180$ such equations.
The non-edges $\overline{\al \be}$ come in two types:
in $100$ cases the line $\overline{\al \be}$ does not lie in a face
of the original dodecahedron, leading to $y_\al y_\be=0$;
if it lies in such a face we get a Del Pezzo equation (in $5\times12=60$
cases).

We summarise:
$$
\matrix{
{\rm type}& {\rm equation} & {\rm conditions}& \# \cr
(1) & x_ix_j  & i\neq j  &  66 \cr
(2) & x_iy_\al & \overline{i\al}{\rm \ not\ an\ edge} & 180\cr
(3) & y_\al y_\be & \overline{\al \be}{\rm \ not\ in\ face} & 100\cr
(4) & y_\al y_\ga-x_i y_\be -x_i^2& {\rm Del\ Pezzo} & 60}
$$

The relations follow from the relations between the generators of 
the Stanley-Reisner ideal of the great stellated dodecahedron, which have 
a particularly simple form: we get a relation
for each pair of equations which have one variable in common.

This gives the following list of relations, where we suppress the
conditions on the indices; they can be deduced from the list of equations.
$$
\matrix{
{\rm type}& {\rm relation} & \# \cr
(1-1)  &  (x_ix_j)x_k-(x_ix_k)x_j     & 440  \cr
(2-1)  &  (x_iy_\al)x_j-(x_ix_j)y_\al   &  1980\cr
(2-2)  &  (x_iy_\al)y_\be - (x_iy_\be)y_\al & 1260 \cr
(3-2)  &  (y_\al y_\be)x_i-(x_iy_\al)y_\be & 1200 \cr
(3-3)  &  (y_\al y_\be)y_\ga -(y_\al y_\ga)y_\be  & 780 \cr
(4-2)  &  (y_\al y_\ga - x_iy_\be -x_i^2)x_j
            -(x_jy_\al)y_\ga+(x_ix_j)(x_i+y_\be)  & 660\cr
(4-3)  &  (y_\al y_\ga - x_iy_\be -x_i^2)y_\de
            -(y_\al y_\de)y_\ga+(x_iy_\de)(x_i+y_\be)  & 860\cr
(4-3)  &  (y_\al y_\ga - x_iy_\be -x_i^2)y_\de
            -(y_\al y_\de-x_jy_\be-x_j^2)y_\ga \hfill\cr
 & \hfill{}+(x_iy_\de)(x_i+y_\be)
            -(x_jy_\ga)(x_j+y_\be)  & 40\cr
(4-4)  &  {\rm relations\ from \ matrices}   & 60
}
$$

We use the equations and relations to compute infinitesimal deformations.
The computations are similar to the case of the tetrahedron
of degree 12. 
To illustrate our methods we prove that the
dodecahedron $X$ has no nontrivial extensions.
This statement means that $X$ is only  a hyperplane section
of the projective cone over it. To prove this we have to show that
the affine cone
$C(X)$ has no deformations of negative degree.

\proclaim Proposition.
$T^1_{C(X)}(-\nu)=0$ for $\nu>0$.

\roep Proof.
As all quadratic equations occur in linear relations we cannot perturb the
equations with constants, so $T^1(-2)=0$.

Now we consider deformations of degree $-1$. We start with equations of
type $(1)$, which we perturb as follows:
$$
x_ix_j+\sum a_{ij}^m x_m + \sum b_{ij}^\al y_\al \;.
$$
The relations $(1-1)$ together with the equations give
$$
a_{ij}^k x_k^2 + \sum_{\overline{k\al}{\rm\ edge}} b_{ij}^\al x_ky_\al =
a_{ik}^j x_j^2 + \sum_{\overline{j\al}{\rm\ edge}} b_{ik}^\al x_jy_\al\;.
$$
This shows that $a_{ij}^k=0$ for $k\notin\{i,j\}$. For each $\al$ we can 
find a $k\notin\{i,j\}$ such that $\overline{k\al}$ is an  edge, so
$b_{ij}^\al=0$ and the deformation has the form
$$
x_ix_j+ a_{ij}^i x_i + a_{ij}^j x_j \;.
$$
We now perturb equations of type $(2)$:
$$
x_iy_\al+\sum a_{i\al}^m x_m + \sum b_{i\al}^\be y_\be 
$$
and use the relations of type $(2-1)$:
$$
a_{i\al}^j x_j^2 + \sum_{\overline{j\be}{\rm\ edge}} b_{i\al}^\be x_j y_\be=
a_{ij}^j x_jy_\al 
$$
to conclude that $a_{i\al}^j=0$ for all $j\neq i$, 
$b_{i\al}^\be=0$ for all $\be\neq\al$ and $b_{i\al}^\al=a_{ij}^j$
for all $j$ such that $\overline{j\al}$ is an edge.
It follows that $a_{ij}^j=a_{ik}^k$ for all $j$ and $k$. Using the
coordinate transformation $\partial_{x_i}$ we may therefore assume
that the equations of type $(1)$ are not perturbed at all, while those
of type $(2)$ have the form $x_iy_\al+a_{i\al}^i x_i$.

Perturbing equations of type $(3)$ in a similar manner as
$y_\al y_\be+\sum a_{\al\be}^m x_m + \sum b_{\al\be}^\ga y_\ga$
we find from the relations $(3-2)$ that
$$
a_{\al\be}^i x_i^2 + \sum_{\overline{i\ga}{\rm\ edge}} b_{\al\be}^\ga x_iy_\ga
=a_{i\al}^i x_iy_\be \;.
$$
So $a_{\al\be}^i=0$, $b_{\al\be}^\ga=0$ for $\ga\notin\{\al,\be\}$ and
$b_{\al\be}^\be=a_{i\al}^i$. The coordinate transformation 
$\partial_{y_\al}$ can be used to eliminate the perturbation
of the equations of type $(2)$. Then also those of type $(3)$ are
not perturbed.

Finally we look at the Del Pezzo equations $(4)$. From the relations 
$(4-2)$ we conclude as before that the only possible perturbations have the
form
$$
y_\al y_\ga - x_iy_\be -x_i^2+a_{\al\ga}^ix_i\;.
$$
Here $\overline{i\al}$, $\overline{i\be}$ and $\overline{i\ga}$
all are edges. This means that we can look at each Del Pezzo 
separately. From the matrix of relations we obtain that 
$a_{\al\ga}^i=0$.
\qed

\proclaim Proposition.
The dodecahedron $X$ with icosahedral symmetry
is $d$-semi\-stable. The space of locally trivial embedded
deformations has dimension $29$.

\roep Proof.
We describe a global section of the sheaf $\ct_X^1$, without proof.
To formulate the result  we use 
the alternative notation  $y_{ijk}$ for $y_\al$, if 
$\overline{i\al}$, $\overline{j\al}$ and $\overline{k\al}$
are the edges involving $\al$.

The equations of type $(1)$ are not deformed, unless $\overline{ij}$
is an edge, in which case  we have
$x_ix_j + d y_{ijp} y_{ijq}$. 
An equation of type $(2)$ is perturbed
to $x_iy_\al+d(x_jy_{ijp}+ x_jy_{ijq}+x_j^2)$
if $\al=(jkl)$ is opposite the edge $[pq]$;
if $ijp$, $ijq$ and $\al=jqr$ are three consecutive vertices, as are
$\al=jqr$, $ijq$ and $iqs$, we get 
$x_iy_\al+d(x_jy_{ijq}+ x_qy_{ijq}+y_{ijq}^2)$
and in all  other cases the equation is not deformed.
$$
\setlength{\unitlength}{.8cm}
\begin{picture}(6,6.5)(-2.5,-3.25)
\put(-1.5,0){\line(1,0){3}}
\put(-1.5,0){\line(-1,2){1}}
\put(-1.5,0){\line(-1,-2){1}}
\put(1.5,0){\line(1,2){1}}
\put(1.5,0){\line(1,-2){1}}
\put(2.5,2){\line(-2,1){2.5}}
\put(2.5,-2){\line(-2,-1){2.5}}
\put(-2.5,2){\line(2,1){2.5}}
\put(-2.5,-2){\line(2,-1){2.5}}
\put(0,1.5){\makebox(0,0){$i$}}
\put(0,-1.5){\makebox(0,0){$j$}}
\put(2.5,0){\makebox(0,0){$p$}}
\put(-2.5,0){\makebox(0,0){$q$}}
\put(0,-3){\makebox(0,0){$\al$}}
\end{picture}
\qquad
\begin{picture}(7,6.5)(-4,-6.5)
\put(-1.5,0){\line(1,0){3}}
\put(-1.5,0){\line(-1,-2){1}}
\put(1.5,0){\line(1,-2){1}}
\put(2.5,-2){\line(-2,-1){2.5}}
\put(-2.5,-2){\line(2,-1){2.5}}
\put(2.5,-2){\line(2,-1){.5}}
\put(-2.5,-2){\line(-2,-1){.5}}
\put(0,-3.25){\line(0,-1){3}}
\put(0,-1.5){\makebox(0,0){$i$}}
\put(1.7,-4.5){\makebox(0,0){$j$}}
\put(-1.7,-4.5){\makebox(0,0){$q$}}
\put(0,-6.5){\makebox(0,0){$\al$}}
\put(2,0){\makebox(0,0){$\be$}}
\put(2.7,-1.5){\makebox(0,0){$p$}}
\put(-2.7,-1.5){\makebox(0,0){$s$}}
\end{picture}
$$
Perturbation of the equations $(3)$
all vanish except when $\al$ and $\be$ are nearest possible:
one can reach $\be$ from $\al$ by passing three edges.
Suppose that the vertices on this path are $\be=ipo$, $ijp$, $ijq$, 
$jqr=\al$. Then we have 
$y_\al y_\be+d(x_iy_{ijq}+x_jy_{ijp}+y_{ijp}y_{ijq})$.
The Del Pezzo equations are not deformed.

We compute locally
near a triple point and look at the chart $y_\al=1$. All variables
can be eliminated except $y_\be$, $y_\ga$ and $y_\de$ such that
$\overline{\al\be}$, $\overline{\al\ga}$ and $\overline{\al\de}$ are
edges, and $x_i$, $x_j$ and $x_k$, where $\al$, $\be$ and $\ga$ lie
on the face $k$, etc. We have nine equations left, of three types:
$x_ix_j+dy_\de$, $x_ky_\de+d(1+x_i+x_j)$ and the Del Pezzo equation
$y_\be y_\ga-x_k-x_k^2$. The last one shows that even $x_k$ can be 
eliminated, as the double curve lies in $x_k=0$. By multiplying the
Del Pezzo equation with $y_\de$ and using the other equations we get 
$$
y_\be y_\ga y_\de +d(1 +x_i +x_j +x_k) \;.
$$
This show that our $d$-deformation indeed represents the class
$[1]\in  H^0(\sier D)=H^0(\ct^1_X)$. 
\qed
 
Prop.~(\ref{ingebed}) gives 
the dimension of the space of locally trivial deformations, 
but it can also computed directly. For each
Del Pezzo we have 5 deformations by multiplying the $x_i$ in a given
column of the defining matrix with a unit of the form $1+\ep_{ij}$.
These deformations are trivial and can also be obtained by multiplying
the $y_\al$ by suitable factors. In total we have 60 such deformations,
but globally we have only 31 diagonal coordinate transformations (32
variables, but we have to subtract one for the Euler vector field).

\proclaim Theorem.
There exists a semistable degeneration of $K3$-surfaces of degree 60
with icosahedral symmetry, whose special fibre is our dodecahedron.

\neu
The rotation group $G_{60}\cong A_5$ of the icosahedron  acts 
symplectically on the general fibre ${\cal X}_t$ and the quotient
${\cal X}_t/G_{60}$
is again a $K3$-surface, with 2 $A_4$, 3 $A_2$ and 4 $A_1$ singularities
\cite{X}. The locus of such surfaces has dimension two in moduli, so
together with a polarisation there is only a curve of such surfaces.
It would be interesting to know this curve. A deformation computation
as above only gives a parametrisation with power series; anyway, the
computation is too complicated.

We can take the quotient of the special fibre, which is our
dodecahedron. Invariants for the icosahedral reflection group are
$$
\eqalign{
X & = \sum x_i \cr
Y_1 & = \sum y_\al \cr
Y_2 & = \sum_{\overline{\al\be}  \ {\rm edge}} y_\al y_\be\;,}
$$
and a skew invariant $Z$ is obtained by taking the $G_{60}$ orbit
of $x_iy_\al y_\be(y_\al - y_\be)$.
After a coordinate transformation $X\mapsto \frac15 X$,
$Y_2\mapsto Y_2+\frac 15 XY_1 + \frac15 X^2$ the quotient is given 
by the equation
$$
\displaylines{\quad
Z^2= - X^2\left(5Y_2^2(4Y_2+8X^2+12XY_1-Y_1^2)
\right.\hfill\cr\hfill
\left.{}+(30Y_1+20X)Y_2X(X^2+XY_1-Y_1^2)) 
+(3X^2+4XY_1)(X^2+XY_1-Y_1^2)^2\right).\;}
$$ 
This is a surface of degree 8 in the weighted projective space
$\P(1,1,2,4)$. These numbers are in Reid's list of famous 95 and the
general $X_8 \subset \P(1,1,2,4)$ is a $K3$-surface with 2 $A_1$
singularities. Our surface has a double line and 
two $A_4$ singularities, at $Y_2=X^2+XY_1-Y_1^2=0$.

\neu Finite groups acting symplectically on $K3$-surfaces
have been classified by Mukai \cite{Mu}, see also \cite{X}. Mukai
gives an example of a $K3$-surface with even an action of the symmetric
group $S_5$: 
$$
\eqalign{
x_1+x_2+x_3+x_4+x_5&=0\cr
x_0^2+x_1^2+x_2^2+x_3^2+x_4^2+x_5^2&=0\cr
x_1^3+x_2^3+x_3^3+x_4^3+x_5^3&=0\;.
}
$$


\bfneu{The Stanley-Reisner ring of the icosahedron.}
A different family of $K3$-surfaces with icosahedral symmetry 
is obtained by smoothing the Stanley-Reisner ring of the icosahedron.
The infinitesimal deformations can be found from \cite[Sect.~4]{AC}
or computed directly with the methods above.
All deformations are unobstructed ($T^2_X=0$). We have 
$T^1_X(\nu)=0$ for $\nu<0$ and $\dim T^1_X(0)=30$.
Furthermore the dimension of $H^0(\Theta_X)$ equals 11,
which fits with the fact that $X$ deforms to smooth $K3$-surfaces
($30-11=19$).

\begin{figure}[h]
\centering
\includegraphics[width=7cm]{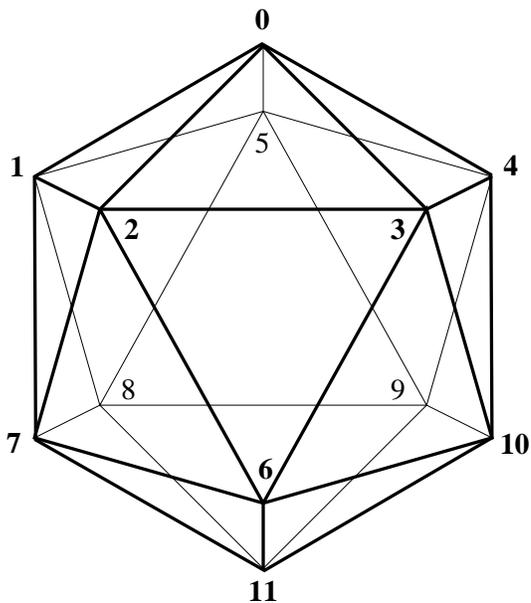}
\caption{Icosahedron.}\label{plaatje}
\end{figure}

We number the vertices as in Fig.~\ref{plaatje}.
Then we have two types of equations, depending on the
distance between vertices.
The infinitesimal deformations are
$$
x_0x_6+\ep_{06}x_2x_3, \qquad x_0x_{11}
$$
By taking all $\ep_{ij}$ equal we get an icosahedral invariant
deformation. The lift to a one-parameter deformation seems to
involve power series of the deformation variable (I computed up to order 
7). Anyway, equations for a $K3$ of degree $20$ are not very illuminating.

As before, this deformation is not semistable, because the total space has
singularities.  Each vertex of the icosahedron gives a singularity,
which is the cone over a pentagon. It is smoothed negatively,
with total space the cone over a Del Pezzo of degree 5.
We resolve these singularities by blowing up.
We introduce 12 Del Pezzo surfaces. The sides are blown up in three points,
giving hexagons. The dual graph of the central fibre is now
a stellated dodecahedron. The object itself 
consists of pentagons and hexagons.  It contains a real homology class,
as described in (\ref{klasse}),
which looks like a football, so our special fibre 
is a complexified football.

\bfneu{A degeneration of degree 12.}
The existence of the two degenerations above with icosahedral symmetry
follows from a deformation argument, but it is too complicated
to give explicit equations. In the football case pentagons arise
because of the singularities of the total space. This suggests
that one can get a degeneration of low degree by blowing down
components of the special fibre. Blowing down means removing
vertices from the dual graph.

We start from the icosahedron (Fig.~\ref{plaatje}) and remove
non-adjacent vertices, say those numbered $0$, $7$ and $10$.
This means breaking the symmetry. The resulting dual graph 
\begin{figure}[h]
\centering
\includegraphics[width=5cm]{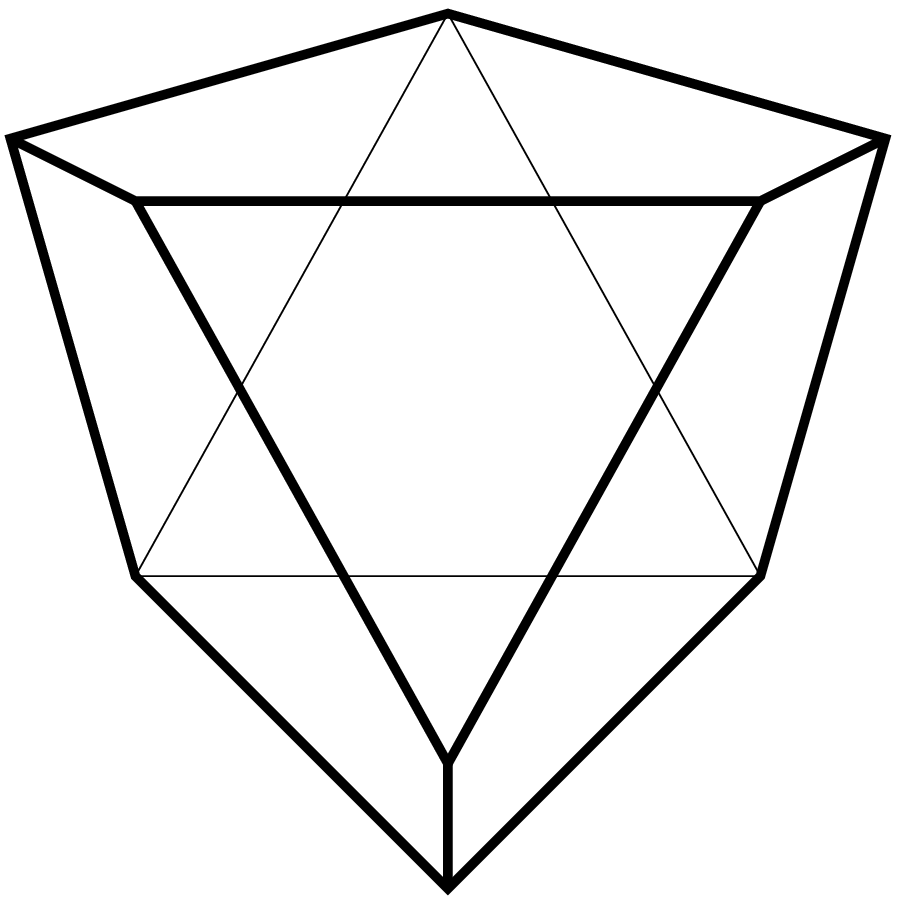}
$\quad$
\includegraphics[width=5cm]{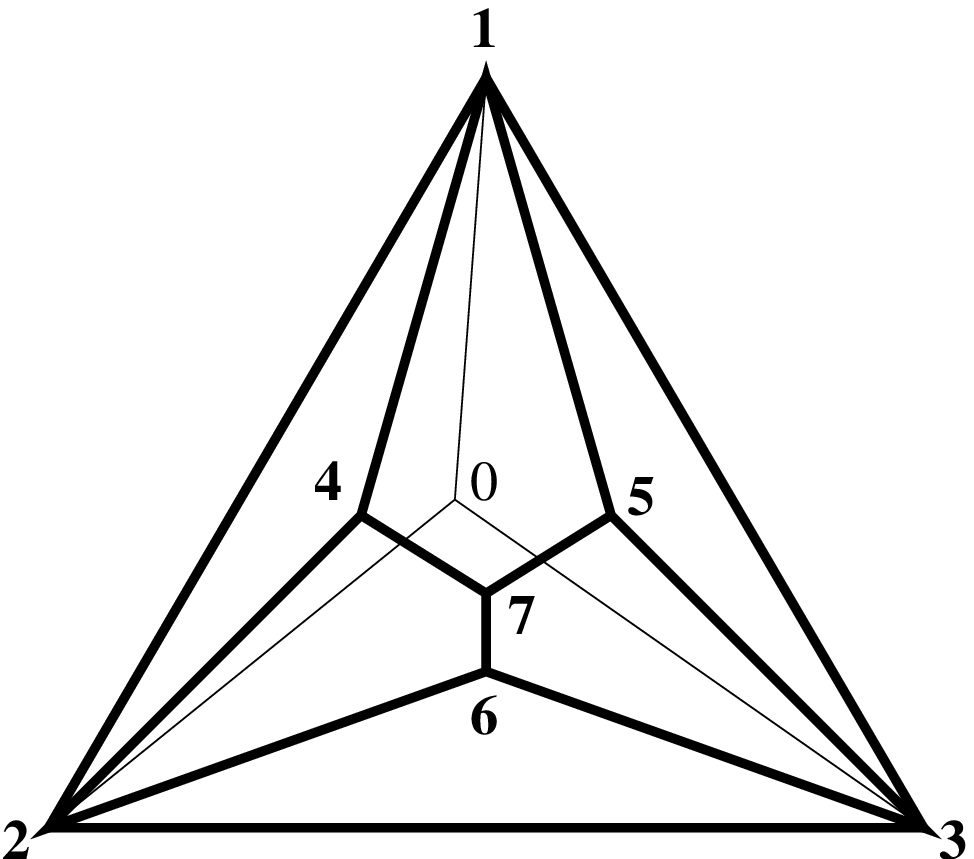}
\caption{Dual graph and its realisation for $X$ of degree 12.}
\label{plaatje2}
\end{figure}
is shown in Fig.~\ref{plaatje2}. Of the double curves on the components
six are triangles and four are rectangles. We realise them on planes
resp.~quadric surfaces. The picture also shows a realisation (as a real
polyhedron). We cannot take the Stanley-Reisner ideal, as the
realisation contains rectangles. 
For those  we take an equation of the form $xy-zt$.
One may think that a rectangle can be triangulated in two ways,
each giving a monomial, which are forced  to be equal.
The result is a surface $X\subset \P^7$ of degree 12. 
With the numbering in the figure we get  the $S_3$-invariant ideal
$$
\displaylines{
x_0x_7, \cr
x_0x_4, \quad x_0x_5, \quad x_0x_6, \cr
x_1x_6, \quad x_2x_6, \quad x_3x_4, \cr
x_1x_7-x_4x_5, \quad x_2x_7-x_4x_6, \quad x_3x_7-x_5x_6, \cr
x_1x_2x_3\;.
}
$$
The next thing to do is to compute the $T^1$ and $T^2$ for the
affine cone $C(X)$ over $X$. This is conveniently done with a
computer algebra program. A computation with {\sl Macaulay\/} \cite{mac}
gives the following result:

\proclaim Lemma.
As $\sier{C(X)}$-module $T^1_{C(X)}$ is generated by eight elements,
represented by the following perturbations of the equations: 
$$
\displaylines{
x_0x_7, \cr
x_0x_4-c_3x_1x_2, \quad x_0x_5-c_2x_1x_3, \quad x_0x_6-c_1x_2x_3, \cr
x_1x_6+b_0x_7+b_1x_6+b_2x_5+b_3x_6, \cr 
x_2x_5+b_0x_7+b_1x_6+b_2x_5+b_3x_6, \cr 
x_3x_4+b_0x_7+b_1x_6+b_2x_5+b_3x_6, \cr
x_1x_7-x_4x_5, \quad x_2x_7-x_4x_6, \quad x_3x_7-x_5x_6, \cr
x_1x_2x_3+ax_0+b_1x_2x_3+b_2x_1x_3+b_3x_1x_2,
}
$$
and $\dim T^2_{C(X)}=2$, concentrated in degree $-2$. 

The quadratic obstruction is given by $a(c_1-c_2)=a(c_1-c_3)=0$.
We conclude that the degree zero deformations are unobstructed.
The base space for $C(X)$ in non-positive degrees has two components.
As we are mainly interested in $S_3$-invariant deformations we consider
only the component with $c_1=c_2=c_3(=:c)$. The component will be obtained
by substituting polynomials for the deformation variables $a$, $b_i$ and
$c$. A computation gives the equations
$$
\displaylines{
x_0x_7+c(b_0x_7+b_1x_6+b_2x_5+b_3x_6+ac), \cr
x_0x_4-cx_1x_2, \quad x_0x_5-cx_1x_3, \quad x_0x_6-cx_2x_3, \cr
x_1x_6+b_0x_7+b_1x_6+b_2x_5+b_3x_6+ac, \cr 
x_2x_5+b_0x_7+b_1x_6+b_2x_5+b_3x_6+ac, \cr 
x_3x_4+b_0x_7+b_1x_6+b_2x_5+b_3x_6+ac, \cr
x_1x_7-x_4x_5, \quad x_2x_7-x_4x_6, \quad x_3x_7-x_5x_6, \cr
x_1x_2x_3+ax_0+b_1x_2x_3+b_2x_1x_3+b_3x_1x_2-
        b_0(b_0x_7+b_1x_6+b_2x_5+b_3x_6+ac),
}
$$
For $c\neq0$ we derive the three equations $x_0x_7-cx_ix_{7-i}$, which show
that we have a hypersurface in the cone over $\P^1\times \P^1\times \P^1$.

\proclaim Proposition. \label{icodual}%
A general one-parameter deformation in degree 0 on the
component described above has a minimal model in $(-1)$-form
with the icosahedron as dual graph for the central fibre.
In particular, this holds for 
$$
\eqalign{
a&=c(x_0^2+x_0(x_1+x_2+x_3)+x_1^2+x_2^2+x_3^2),\cr
b_0&=cx_7,\cr 
b_1&=c(x_2+x_3+x_6+x_7), \cr 
b_2&=c(x_1+x_3+x_5+x_7), \cr 
b_3&=c(x_1+x_2+x_4+x_7).
}
$$ 

\roep Proof.
One first checks that the general fibre is a smooth $K3$-surface.
For this it suffices to look at the hypersurface in 
$\P^1\times \P^1\times \P^1$.

In the particular example the total space has at the origin of the
affine chart $x_1=1$ a singularity, which is isomorphic to the cone
over the Del Pezzo surface of degree 5, as it should be: the point
to check is that we indeed have a generic local deformation.
Furthermore there are 18 singularities of type $A_1$.
On the $(x_1,x_4)$-line we have the point $x_1+x_4=0$.
On the $(x_0,x_1)$-line we have two points, given by $x_0^2+x_0x_1+x_1^2$,
and on  the $(x_7,x_7)$-line the two points  $x_6^2+x_6x_7+x_7^2$.
The other singular points are found by symmetry.

By blowing up the three singularities of multiplicity 5 and 
making a small resolution of the $A_1$-points we get a smooth total
space. To obtain the $(-1)$-form one has to place one exceptional curve
on either component in case the double line contains two singularities.
If there is only one, the exceptional curve should lie on the
triangle component.
\qed

\parskip=0pt plus 1pt

\vfill

{\obeylines
\parindent=0pt
Address of the author:
Matematiska institutionen
G\"oteborgs universitet
Chalmers tekniska h\"ogskola
SE 412 96 G\"oteborg, Sweden
email: stevens@math.chalmers.se}

\end{document}